\input amstex
\mag=\magstep1
\documentstyle{amsppt}
\nologo
\define\BaseLineSkip{\baselineskip=13.5pt}  
\NoBlackBoxes
\rightheadtext\nofrills{\eightpoint Determinant expressions for hyperelliptic functions}
\leftheadtext\nofrills{\eightpoint Yoshihiro \^Onishi}
\pagewidth{138mm} 
\hcorrection{-0pt} 
\vcorrection{-35pt} 
\define\inbox#1{$\boxed{\text{#1}}$}
\def\fp{\flushpar}
\define\dint{\dsize\int}
\define\underbarl#1{\lower 1.4pt \hbox{\underbar{\raise 1.4pt \hbox{#1}}}}

\define\tp#1{\negthinspace\ ^t#1}
\define\w{\omega}
\define\lr#1{^{\sssize\left(#1\right)}}
\define\br#1{^{\sssize\left[#1\right]}}
\define\do#1{^{\sssize\left<#1\right>}}
\define\lqq{\lq\lq}
\define\C#1#2{{{#1}\choose{#2}}}

\define\dg#1{(d^{\circ}\geqq#1)}
\define\Dg#1#2{(d^{\circ}(#1)\geqq#2)}
\define\Dxd{\frac{d}{dx}}
\define\Dx#1{\left(\frac{d}{dx}\right)^{#1}}
\def\rank{\text{\rm rank\,}}
\define\nullhbox{\hbox{\vrule height0pt depth 0pt width 1pt}}
\define\qedright{\null\nobreak\leaders\nullhbox\hskip10pt plus1filll\ \qed\par}


\font\tenptmbit=cmmib10 
\define\bk#1{\text{\tenptmbit{#1}}}

\font\elevenptrm=cmr10 scaled \magstephalf

\font\Large=cmr10 scaled \magstep5
\topmatter
\title  \nofrills\elevenptrm 
DETERMINANT EXPRESSIONS FOR HYPERELLIPTIC FUNCTIONS \\
\elevenptrm (with an Appendix by Shigeki Matsutani)
\endtitle
\author {YOSHIHIRO \^ONISHI} \\
\endauthor
\abstract\nofrills{\it Abstruct}\ \ 
In this paper we give quite pretty generalization of 
the formula of Frobenius-Stickelberger to all hyperelliptic curves.  
The formula of Kiepert type is also obtained by limiting process 
from this generalization. 
In Appendix a determinant expression of D.G. Cantor is also given. 
\endabstract
\endtopmatter
\document

\TagsOnRight
\document

\fp
\baselineskip=10pt 

\vskip 3pt
\fp
{\eightpoint 2000 \it Mathematics subject classification: {\rm 11G30, 11G10, 14H45}}
\fp

\BaseLineSkip

\subheading{Introduction}

\vskip 5pt

\fp
Let  $\sigma(u)$  and  $\wp(u)$  be the usual functions 
in the theory of elliptic functions.  
The following two formulae were found in the nineteenth-century.  
First one is 
  $$
  \aligned
  (-1)^{(n-1)(n-2)/2}1!2!&\cdots (n-1)!
   \frac{\sigma(u\lr{1}+u\lr{2}+\cdots+u\lr{n})\prod_{i<j}\sigma(u\lr{i}-u\lr{j})}
        {\sigma(u\lr{1})^n\sigma(u\lr{2})^n\cdots\sigma(u\lr{n})^n} \\
  &= \left|\matrix
      1 &  \wp(u\lr{1}) & \wp'(u\lr{1}) & \wp''(u\lr{1}) & \cdots & \wp\lr{n-2}(u\lr{1}) \\
      1 &  \wp(u\lr{2}) & \wp'(u\lr{2}) & \wp''(u\lr{2}) & \cdots & \wp\lr{n-2}(u\lr{2}) \\
 \vdots &  \vdots   & \vdots    & \vdots     & \ddots & \vdots           \\
      1 &  \wp(u\lr{n}) & \wp'(u\lr{n}) & \wp''(u\lr{n}) & \cdots & \wp\lr{n-2}(u\lr{n}) \\
   \endmatrix\right|.
  \endaligned
  \tag 0.1
  $$
This formula appeared in  
the paper of Frobenius and Stickelberger \cite{{\bf 11}}.  
Second one is 
  $$
  (-1)^{n-1}(1!2!\cdots (n-1)!)^2
   \frac{\sigma(nu)}{\sigma(u)^{n^2}}
  = \left|\matrix
    \wp'        & \wp''     &  \cdots  & \wp\lr{n-1}  \\
    \wp''       & \wp'''    &  \cdots  & \wp\lr{n}    \\
    \vdots      & \vdots    &  \ddots  & \vdots       \\
    \wp\lr{n-1} & \wp\lr{n} &  \cdots  & \wp\lr{2n-3} \\
   \endmatrix\right|(u).   
  \tag 0.2
  $$
Although this formula can be obtained by a limiting process from (0.1), 
it was found before \cite{{\bf 11}} by the paper of Kiepert \cite{{\bf 13}}.  

If we set  $y(u)=\frac 12 \wp'(u)$  and  $x(u)=\wp(u)$,   
then we have an equation  $y(u)^2=x(u)^3+\cdots$, that is a defining equation 
of the elliptic curve 
to which the functions  $\wp(u)$  and  $\sigma(u)$  are attached.  
Here the complex number  $u$  and the coordinate  $(x(u), y(u))$  
correspond by the equality  
  $$
  u=\dint_{\infty}^{(x(u), y(u))}\dfrac{dx}{2y}.
  $$
Then (0.1) and (0.2) is easily rewritten as
  $$
  \aligned
  (-1)^{(n-1)(n-2)/2}
   &\frac{\sigma(u\lr{1}+u\lr{2}+\cdots+u\lr{n})\prod_{i<j}\sigma(u\lr{i}-u\lr{j})}
        {\sigma(u\lr{1})^n\sigma(u\lr{2})^n\cdots\sigma(u\lr{n})^n} \\
  &=\left|\matrix
     1 &  x(u\lr{1}) & y(u\lr{1}) & x^2(u\lr{1}) & yx(u\lr{1}) & x^3(u\lr{1}) & \cdots  \\
     1 &  x(u\lr{2}) & y(u\lr{2}) & x^2(u\lr{2}) & yx(u\lr{2}) & x^3(u\lr{2}) & \cdots  \\
\vdots &  \vdots & \vdots   & \vdots & \vdots   & \vdots  & \ddots  \\
     1 &  x(u\lr{n}) & y(u\lr{n}) & x^2(u\lr{n}) & yx(u\lr{n}) & x^3(u\lr{n}) & \cdots  \\
   \endmatrix\right|
  \endaligned
  \tag 0.3
  $$
and
  $$
  \align
  (-1)^{n-1}1!2!&\cdots (n-1)!
  \frac{\sigma(nu)}{\sigma(u)^{n^2}} \\
  &=\left|\matrix
    x'           & y'            & (x^2)'  
  & (yx)'        & (x^3)'        & \cdots  \\
    x''          & y''           & (x^2)'' 
  & (yx)''       & (x^3)''       & \cdots  \\
    \vdots       & \vdots        & \vdots   
  & \vdots       & \vdots        & \ddots  \\
    x\lr{n-1}    & y\lr{n-1}     & (x^2)\lr{n-1} 
  & (yx)\lr{n-1} & (x^3)\lr{n-1} & \cdots  \\
  \endmatrix\right|(u), 
  \tag 0.4
  \endalign
  $$
respectively.  

The author recently gave a generalization of 
the formulae (0.3) and (0.4) to the case of genus two 
in \cite{{\bf 18}} and to the case of genus three in \cite{{\bf 19}}. 
The aim of this paper is to give 
a quite natural genaralization of  (0.3),  (0.4)  
and the results in \cite{{\bf 18}},  \cite{{\bf 19}}  
to all of the hyperelliptic curves  (see Theorem 7.2 and Theorem 8.3).  
The idea of our generalization arises from 
the unique paper \cite{{\bf 12}} of D. Grant 
and it can be said in a phrase, 
\lqq{\it Think not on the Jacobian but on the curve itself.}"  

Fay's famous formula (44) in p.33 of \cite{{\bf 10}} 
is a formula on the Jacobian, and is regarded 
as a generalization of (0.3).  
The author do not know whether his formula is able to yield 
a generalization of (0.4).  
Our formula is quite pretty in a comparison with Fay's one 
and naturally gives a generalization of (0.4).  
To connect Fay's formula with ours would be a problem 
\footnote{Recently the paper \cite{{\bf 9}} appeared. 
This paper seems to investigate this problem} in future. 

Now we prepare the minimal fundamentals to explain our results.  
Let  $f(x)$  be a monic polynomial of  $x$  of degree  $2g+1$  
with  $g$  a positive integer.  
Assume that  $f(x)=0$  has no multiple roots.   
Let  $C$  be the hyperelliptic curve defined by  $y^2=f(x)$.  
Then  $C$  is of genus  $g$  and it is ramified at infinity.  
We denote by $\infty$  the unique point at infinity.  
Let  $\bold C^g$  be the coordinate space of all vectors 
of integrals 
  $$
  \left(\int_{\infty}^{P_1}+\cdots+\int_{\infty}^{P_g}\right)
  (\frac1{2y}, \ \frac{x}{2y}, \ \ldots,  \ \frac{x^{g-1}}{2y})dx
  $$  
of the first kind for  $P_j\in C$.  
Let  $\Lambda\subset \bold C^g$  be the lattice of their periods.  
So  $\bold C^g/\Lambda$  is the Jacobian variety of  $C$.  
We denote the canonical map by  $\kappa: \bold C^g \rightarrow \bold C^g/\Lambda$.  
We have an embedding  $\iota:C\hookrightarrow \bold C^g/\Lambda$  
defined 
by  $P\mapsto (\int_{\infty}^P\frac{dx}{2y}, \ 
               \int_{\infty}^P\frac{xdx}{2y}, \  
               \ldots, \ 
               \int_{\infty}^P\frac{x^{g-1}dx}{2y})$ mod $\Lambda$. 
Therefore  $\iota(\infty)=(0, 0, \ldots, 0)\in \bold C^g/\Lambda$.  
We regard an algebraic function on  $C$,  
that we call a {\it hyperelliptic function} in this article,  
as a function on a universal Abelian covering  $\kappa^{-1}\iota(C)$  
($\subset \bold C^g$) of  $C$.  If  $u=(u_1, \ \ldots, \ u_g)$  
is in  $\kappa^{-1}\iota(C)$,  we denote by  $(x(u), y(u))$  
the coordinate of the corresponding point on  $C$  by 
  $$
  u=
  \int_{\infty}^{(x(u), y(u))}
  (\frac{1}{2y}, \ \frac{x}{2y}, \ \ldots, \ \frac{x^{g-1}}{2y})dx 
  $$
with appropriate choice of a path of the integrals.  
Needless to say, we have  $(x(0, 0, \ldots, 0), \allowmathbreak y(0, 0, \ldots, 0))=\infty$.   

Our new standing point of view is characterized by 
the following three comprehension about the formulae (0.3) and (0.4).  
Firstly, the sequence of functions of  $u$  
whose values at  $u=u_j$  are displayed in the  $(j+1)^{\text{st}}$  row 
of the determinant of (0.3)  is just a sequence of the monomials 
of  $x(u)$  and  $y(u)$  displayed according to the order 
of their poles at  $u=0$.  
Secondly, the two sides of (0.3) as a function of  $u=u_0$  
and the those of (0.3) should be regarded as functions 
defined on the universal (Abelian) covering space  $\bold C$  
{\it not} of the Jacobian variety but {\it of the elliptic curve}.  
Thirdly, the expression of the left hand side of (0.4) states 
the function of the two sides themselves of (0.4) is characterized 
as an elliptic function such that its zeroes are exactly the points 
different from  $\infty$  
whose $n$-plication is just on the standard theta divisor 
in the Jacobian of the curve, and such that its pole is only at  $\infty$.  
In the case of the elliptic curve above, 
the standard theta divisor is just the point at infinity.  

Surprisingly enough, these three comprehensions just invent 
a good generalization for hyperelliptic curves.  
More concreatly, for  $n\geqq g$  
our generalization of (0.4) is obtained by replacing 
the sequence of the right hand side by the sequence
  $$
  1,\  x(u),\  x^2(u),\  \ldots, \ x^g(u),\  y(u),\  x^{g+1}(u),\  yx(u),\  \ldots,
  $$  
where  $u=(u_1, u_2, \ldots, u_g)$  is on  $\kappa^{-1}\iota(C)$, 
of the monomials of  $x(u)$  and  $y(u)$  displayed according to 
the order of their poles at  $u=(0, 0, \ldots, 0)$  with replacing the derivatives 
with respect to  $u\in\bold C$  by 
those with respect to  $u_1$  along  $\kappa^{-1}\iota(C)$;  
and the left hand side of (0.4) by  
  $$
  1!2!\cdots (n-1)!\sigma(nu)/\sigma_{\sharp}(u)^{n^2}, 
  \tag 0.5
  $$  
where  $n\geqq g$, $\sigma(u)=\sigma(u_1, u_2, \ldots, u_g)$  
is a well-tuned Riemann theta series which is a natural generalization of 
elliptic  $\sigma(u)$  and  $\sigma_{\sharp}$  is defined as in the table below:  

$ $

\centerline{
\vbox{\offinterlineskip 
\halign{
\strut# & \ \hfil # \hfil \ &&
      # & \ \hfil # \hfil \ \cr 
& genus $g$            \hfill\vrule &
&       1        &
&       2        &&         3      &&        4        &&       5       &
&       6        &&         7      &&        8        &&    $\cdots$   \cr
\noalign{\hrule}
& \hfill $\sigma_{\sharp}$ \hfill\vrule &
& $\sigma$       &
& $\sigma_2$     && $\sigma_2$     && $\sigma_{24}$   && $\sigma_{24}$ &
& $\sigma_{246}$ && $\sigma_{246}$ && $\sigma_{2468}$ &&    $\cdots$   \cr
}}}

$ $

\fp
where 
  $$
  \sigma_{ij\cdots\ell}(u)
  =
  \tfrac{\partial}{\partial u_i}
  \tfrac{\partial}{\partial u_j}
  \cdots
  \tfrac{\partial}{\partial u_{\ell}}
  \sigma(u).
  $$
The function (0.5) is a good generalization of the $n$-{\it division polynomial} \ 
of an elliptic curve as mentioned in Remark 8.4 below.  
For the case  $n\leqq g$, we need slight modification as in Theorem 7.2(1) and  8.3(1).  
As a function on  $\kappa^{-1}\iota(C)$,  $\sigma_{\sharp}(u)$  has the zeroes at 
the points   $\kappa^{-1}\iota(\infty)$  and no zero elsewhere (Proposition 6.6(1)).  
This property is exactly the same as the elliptic  $\sigma(u)$.  
The hyperelliptic function (0.5) 
can be regarded as a function on  $\bold C^g$  via theta functions.  
Although this function on  $\bold C^g$  is no longer 
a function on the Jacobian, 
it is indeed expressed simply in terms of theta functions and is treated 
really similar way to the elliptic functions.  

The most difficult problem was to find the left hand side of 
the expected our generalization of (0.3).  
For simplicity we assume  $n\geqq g$.  
The answer is remarkably pretty and is 
  $$
  \frac{\sigma(u\lr{1}+u\lr{2}+\cdots+u\lr{n})\prod_{i<j}\sigma_{\flat}(u\lr{i}-u\lr{j})}
  {\sigma_{\sharp}(u\lr{1})^n
   \sigma_{\sharp}(u\lr{2})^n
   \cdots
   \sigma_{\sharp}(u\lr{n})^n}, 
  $$
where  $u\lr{j}=(u\lr{j}_1,  u\lr{j}_2,  \ldots, u\lr{j}_g)$  ($j=1$, $\ldots$, $n$)    
are variables on  $\kappa^{-1}\iota(C)$  
and  $\sigma_{\flat}(u)$  is defined as in the table below:

\vskip 5pt

\centerline{
\vbox{\offinterlineskip 
\halign{
\strut# & \ \hfil # \hfil \ &&
      # & \ \hfil # \hfil \ \cr 
& genus  $g$  \hfill\vrule &
&      1     &
&      2     &&      3     &&        4      &&       5       &&       6        &&         7      &&       8         && $\cdots$ \cr
\noalign{\hrule}
& \hfill $\sigma_{\flat}$  \hfill\vrule &
& $\sigma$   &
& $\sigma$   && $\sigma_3$ && $\sigma_3$    && $\sigma_{35}$ && $\sigma_{35}$  && $\sigma_{357}$ && $\sigma_{357}$  && $\cdots$ \cr
}}}

\vskip 5pt

If we once find this, we can prove the formula, roughly speaking, 
by comparing the divisors of the two sides.  
As the formula (0.4) is obtained by 
a limiting process from (0.3), our generalization of (0.4) is 
obtained by similar limitting process from the generalization of (0.3). 

Cantor \cite{{\bf 8}} gave another determiant expression of 
the function that is characterized in the third comprehension above 
for any hyperelliptic curve.  
The expression of Cantor should be seen as a generalization 
of a formula due to Brioschi (see \cite{{\bf 5}}, p.770, $\ell$.3).  

Concerning the paper \cite{{\bf 18}} 
Matsutani pointed out that (0.4) can be 
generalized to all of hyperelliptic curves 
and he proved such the formula is equivalent to Cantor's one.  
He kindly permitted the author to include his proof as an Appendix 
in this paper.  

Matsutani's observation made 
the author to start working on 
an extension of (0.3) for all hyperelliptic curves.   
The method of this paper is entirely different from 
that of  \cite{{\bf 18}}  and  \cite{{\bf 19}}.  
The method gives probably one of the simplest way 
and is based on the paper cited as \cite{{\bf 7}}.  
At the beginning of this reserch the author computed several cases 
by low blow method as in \cite{{\bf 18}} and \cite{{\bf 19}}.  
When Theorem 7.2 was still a conjecture 
Professor V.Z. Enol'skii suggested the author 
that to prove the conjecture it would be 
important to investigate the leading terms of the sigma function 
as in \cite{{\bf 7}}.  

Now we mention the idea of the proof.  
When the curve  $C$  which is defined by  $y^2=f(x)$  deforms to 
a singular curve  $y^2=x^{2g+1}$  the canonical limit of 
the function  $\sigma(u)$  is known to be so called Schur polynomial in 
the theory of representation of the symmetric group.  
The paper \cite{{\bf 7}} treated this fact quite explicitly. 
Such the limit polynomial is called the {\it Schur-Weierstrass polynomial}~ in that paper. 
For our argument, we need a slight improvement of this fact (see Section 4).  
To prove our formula of Frobenius-Stickelberger type 
by induction on the number of variables  $u\lr{j}$  
we need relations to connect to each factor of the numerator 
to a factor of the denominator in the left hand side of Theorem 7.2.  
So after proving such the connection about the Schur-Weierstrass polynomial 
as explained in Section 2,   
we will lift the connection to the case of the sigma function as 
in Section 6.  
To lift so we need additional facts on vanishing of some derivatives 
of the sigma function as is described in Section 5. 

The results from Section 1 to Section 6 are easily generalized to 
quite wide family of algebraic curves.  
Such curves are called  $(n,s)$-curves in the paper \cite{{\bf 7}}. 
Unfortunitely the standard theta divisor or every standard theta subvariety,  
i.e. an image of the symmetric product of some copies of the curve,  
in the case of such general curve, 
is {\it not symmetric} in the Jacobian.  
Here the word \lqq standard" means that 
the embbeding of the curve to its Jacobian variety 
send the point at infinity to the origin.  
Hence they have no involution and a naive generalization of ours 
ended in failure.  
To go further it would give a hint that we compare our formulae closely 
with Fay's one.  

There are also various generalizations of (0.1) (or (0.3)) 
in the case of genus two different from our line.  
If the reader is interested in them, 
he should be refered to Introduction of \cite{{\bf 18}}.


\vskip 15pt

\subheading{Convention}
\vskip 5pt
\fp
We use the following notations throughout the paper.  
We denote, as usual, by  $\bold Z$, $\bold Q$, $\bold R$  
and  $\bold C$ 
the ring of rational integers, 
the field of rational numbers, 
the field of real numbers 
and the field of complex numbers, respectively.  
In an expression of the Laurent expansion of a function, 
the symbol  $(d^{\circ}(z_1, z_2, \ldots, z_m)\geqq n)$  stands for 
the terms of total degree at least  $n$  with respect to 
the given variables  $z_1$, $z_2$, $\ldots$, $z_m$. 
This notation {\it never}~ means that the terms are monomials 
only of  $z_1$, $\ldots$, $z_m$.  
When the variable or the least total degree are clear from the context, 
we simply denote them by  $(d^{\circ}\geqq n)$  or the dots \lq\lq $\ldots$". 

We will often omit zero entries from a matrix.  
For a simplicity we will ocasionally denote a matrix entry with an asterisk. 

For cross references in this paper, we indicate a formula as (1.2),  
and each of Lemmas, Propositions, Theorems and Remarks also as 1.2.  

\vskip 20pt

\subheading{Contents}
\vskip 5pt
\fp
1. The Schur-Weierstrass polynomial
\fp
2. Derivatives of the Schur-Weierstrass polynomial
\fp
3. Hypelliptic Functions
\fp
4. The Schur-Weierstrass polynomial and the sigma function
\fp
5. Vanishing structure of the sigma functions and of its derivatives
\fp
6. Special derivatives of the sigma function
\fp
7. The Frobenius-Stickelberger type formula
\fp
8. The Kiepert type formula
\fp
Appendix. Connection of the formulae of Cantor-Brioschi and of Kiepert type
\fp $\qquad$          (by S. Matsutani) 

\newpage

\subheading{1. The Schur-Weierstrass Polynomial}

\vskip 5pt
\fp
We begin with a review of fundamentals on Schur-Weierstrass polynomials.  
Our main references are \cite{{\bf 14}} and \cite{{\bf 7}}.  

Let  $g$  be a fixed positive integer.  
Let  $u\lr{1}_g$, $\ldots$, $u\lr{g}_g$  are indeterminates.  
We fix  $n$  ($0\leqq n\leqq g$) and 
we denote by  $\bk{u}_g$  the set of variables  $u\lr{1}_g$, $\ldots$, $u\lr{n}_g$.  
For each  $k\geqq 0$  we denote by  $(-1)^kU\br{n}_k(\bk{u}_g)$ 
the  $k^{\text{th}}$  {\it complete symmetric function}, namely 
the sum of all monomials of total degree  $k$  of 
the variables  $u\lr{1}_g$, $\ldots$, $u\lr{n}_g$.  
We will emphasis by the superscript  $[n]$  
that  $U\br{n}_k(\bk{u}_g)$  is a function 
of a set of  $n$  variables  $\bk{u}_g$.  

We now consider the determinant 
  $$
  |U\br{g}_{g-2i+j+1}(\bk{u}_g)|_{1\leqq i, j \leqq g}.
  $$
If we write simply  $U_k=U\br{g}_k(\bk{u}_g)$, this is concretely of the form 
  $$
  \left|
  \matrix
  U_g     & U_{g+1} & U_{g+2} & \cdots & U_{2g-2} & U_{2g-1} \\
  U_{g-2} & U_{g-1} & U_g     & \cdots & U_{2g-4} & U_{2g-3} \\
  \vdots  & \vdots  & \vdots  & \ddots & \vdots   & \vdots   \\
  U_1     & U_2     & U_3     & \cdots & *        & *        \\
          & U_0     & U_1     & \cdots & *        & *        \\
          &         &         & \ddots & \vdots   & \vdots   \\
          &         &         &        & U_0      & U_1      
  \endmatrix
  \right|
  $$
for odd  $g$, or 
  $$
  \left|
  \matrix
  U_g     & U_{g+1} & U_{g+2} & U_{g+3} & \cdots & U_{2g-2} & U_{2g-1} \\
  U_{g-2} & U_{g-1} & U_g     & U_{g+1} & \cdots & U_{2g-4} & U_{2g-3} \\
  \vdots  & \vdots  & \vdots  & \vdots  & \ddots & \vdots   & \vdots   \\
  U_0     & U_1     & U_2     & U_3     & \cdots & *        & *        \\
          &         & U_0     & U_1     & \cdots & *        & *        \\
          &         &         &         & \ddots & \vdots   & \vdots   \\
          &         &         &         &        & U_0      & U_1      
  \endmatrix
  \right|
  $$
for even  $g$. 

Let 
  $$
  p_j:=\sum_{i=1}^g (u\lr{i}_g)^j. 
  $$
Since the formula \cite{{\bf 14}}, p.29, $\ell$. $-4$  just states 
  $$
  U\br{g}_k(\bk{u}_g)
  =\frac{1}{k!}\left|
  \matrix
  -p_1     &   1       &          &        &      \\
  -p_2     &  -p_1     &   2       &        &      \\
   \vdots  &   \vdots  &   \vdots  & \ddots &      \\
  -p_{k-1} &  -p_{k-2} &  -p_{k-3} & \cdots &   k-1 \\ 
  -p_k     &  -p_{k-1} &  -p_{k-2} & \cdots &  -p_1     
  \endmatrix
  \right|, 
  $$
for  $k=1$, $\ldots$, $g$,   
we see that  $|U\br{g}_k(\bk{u}_g)|$  coincides 
the Schur-Weierstrass polynomial  $S_{2,2g+1}$  
of  \cite{{\bf 7}}, Theorem 4.3.  

In the sequel we denote 
  $$
  \align
  &u\lr{i}_j:=\tfrac1{2(g-j)-1}(u\lr{i}_g)^{2(g-j)-1},  \\ 
  &u\lr{i}:=(u\lr{i}_1, \ldots, u\lr{i}_g), \\
  &u_j:=u\lr{1}_j+u\lr{2}_j+\cdots+u\lr{g}_j, \\
  &u:=u\lr{1}+u\lr{2}+\cdots+u\lr{g}=(u_1, u_2, \ldots, u_g).
  \tag 1.1
  \endalign
  $$
Besides  $|U\br{g}_{g-2i+j+1}(\bk{u}_g)|$  is obviously a polynomial 
of  $u\lr{1}_g$, $\ldots$, $u\lr{g}_g$,  we can prove the following.  

\proclaim{Proposition 1.1} 
The polynomial  $|U\br{g}_{g-2i+j+1}(\bk{u}_g)|_{1\leqq i, j\leqq g}$  above  
is completely determined by the  $g$  values 
of  $u_1$, $u_2$, $\ldots$, $u_g$  defined above.  
\endproclaim

\demo{\it Proof} 
See \cite{{\bf 7}}, p.86, Theorem 4.1.  
\qedright
\enddemo

\fp
So we may write in this paper that
  $$
  S(u):=|U\br{g}_{g-2i+j+1}(\bk{u}_g)|_{1\leqq i, j \leqq g}. 
  $$
This polynomial  $S(u)$  is called the {\it Schur-Weierstrass polynomial}.  

We introduce a weight which is defined by taking  
the weight of  $u_j$  being  $2(g-j)+1$. 
We call this {\it Sato weight}.  
It easy to see that  $S(u)$  is of homogeneous with respect to Sato weight, 
and has Sato weight  $\frac12g(g+1)$.

Let  $m$  be a fixed positive integer 
and  $\xi_1$,  $\ldots$, $\xi_m$  be indeterminates.  
We denote  $(-1)^k$  times the sum of monomials of  
total degree  $k$  of  $\xi_1$, $\ldots$, $\xi_m$   by  $\rho_k(\xi)$,  
where  $\xi$  means the set of  $\xi_1$, $\ldots$, $\xi_m$. 

\definition{Definition 1.2} 
Let  $m$,  $\xi_i$,  and  $\rho_k(\xi)$  are as above.  
A matrix whose every row is successive  $(m+1)$  terms of the sequence
  $$
  \text{$\ldots$, $0$, $0$, $1$, $\rho_1(\xi)$, $\rho_2(\xi)$, $\ldots$}
  $$
except the set of terms 
  $$
  \text{$0$, $\ldots$, $0$, $0$, $1$}
  $$
is called a {\it fundamental matrix without a simple row} 
with respect to  $\xi_1$, $\cdots$, $\xi_m$.  
\enddefinition

\fp
The following is used several times in Section 2.  

\proclaim{Lemma 1.3} 
Let  $m$,  $\xi_i$,  and  $\rho_k(\xi)$  are as above.  
Let  $M$  be a fundamental matrix without a simple row 
with respect to  $\xi_1$, $\ldots$, $\xi_m$.  
We denote by  $\varepsilon_j(\xi)$, 
the elementary symmetric function of  $\xi_1$, $\ldots$, $\xi_m$  
of degree  $j$.  
Then we have
  $$
  M
  \left[\matrix
  \varepsilon_m(\xi) \\
  \varepsilon_{m-1}(\xi) \\
  \vdots \\
  \varepsilon_1(\xi) \\
  1 \\
  \endmatrix\right]
  =
  \left[\matrix
  0 \\ 0 \\ \vdots \\ 0 \\0 
  \endmatrix\right].
  $$
\endproclaim

\demo{\it Proof} 
See \cite{{\bf 14}}, p.21, ($2.6'$).  
\qedright
\enddemo

Although the following will be not used in this paper explicitly,  
it is deeply related to Riemann singularity theorem which is mentioned in Section 5. 
So we display here.  

\proclaim{Lemma 1.4} 
As a polynomial of  $u\lr{1}_g$, $\ldots$, $u\lr{g-1}_g$  we have
  $$
  S(u\lr{1}+\cdots+u\lr{g-1})=0
  $$
identically. 
\endproclaim

\fp{\it Proof}.  
This formula follows from 1.3 by setting  $m=g-1$  and  $M$  to be  
the matrix whose detrerminant expresses  $S(u\lr{1}+\cdots+u\lr{g-1})$.  
\qedright

\newpage

\subheading{2. Derivatives of the Schur-Weierstrass Polynomial}

\vskip 5pt
\fp
We will discuss with some derivatives of Schur-Weierstrass polynomial in order  
to investigate the coresponding derivatives of the sigma function in Section 6.  

\definition{Definition 2.1}
For an integer  $n$  with   $1\leqq n \leqq g$,   
we denote by  $\natural^n$  the set of positive integers  $i$  such 
that  $n+1\leqq i \leqq g$  with  $i\equiv n+1 \mod{2}$.  
\enddefinition

\definition{Definition 2.2}
We denote by  $S_{\natural^n}(u)$  the derivative
  $$
  \left(\prod_{i\in \natural^n}\frac{\partial}{\partial u_i}\right)\negthinspace S(u).  
  $$
Moreover we particularly denote
  $$
  S_{\sharp}(u)=S_{\natural^1}(u), \ \ 
  S_{\flat}(u)=S_{\natural^2}(u).
  $$
\enddefinition

In this Section we denote, as in Section 1, 
  $$
  p_k:=\sum_{i=1}^g (u\lr{i}_g)^k. 
  $$
Hence
  $$
  u_{j}=-\tfrac1{2g-2j+1}p_{2g-2j+1} 
  $$
for  $j=1$,  $\ldots$,  $g$. 

\proclaim{Lemma 2.3}
If we regard a polynomial 
of  $U\br{g}_1(\bk{u}_g)$,  $U\br{g}_2(\bk{u}_g)$, $\ldots$   
as a polynomial of  $p_1$,  $\ldots$,  $p_{2g-1}$,  we have
  $$
  k\frac{\partial}{\partial p_k}
  =(-1)^k\sum_{r\geqq 0}U_r\frac{\partial}{\partial U_{k+r}},
  $$
where we simply write  $U_j=U\br{g}_j(\bk{u}_g)$.  
\endproclaim

\demo{\it Proof}.  
See \cite{{\bf 14}}, p.76. 
\qedright
\enddemo

Now we continue to write  $U_j=U\br{g}_j(\bk{u}_g)$.  
This formula states that   $(-1)^{k-1} k(\partial/\partial p_k)S(u)$  
is the sum of the determinants obtained by \lqq sifting by  $k$ " one of the rows
to the right direction of the matrix of the determinant expression of  $S(u)$.  
If we pay attention to a term of the sum 
which was suffered the sifting by  $k$  at  $i^{\text{th}}$ row, 
we say that  $U_k$  in the  $i^{\text{th}}$ row (or in the  $j^{\text{th}}$  column if no confusion airses) 
{\it is changed to}  $U_0$  by using  $k(\partial/\partial p_k)$.  

\proclaim{Proposition 2.4}
Let  $n$  be an integer such that  $1\leqq n \leqq g-1$.  Let
  $$
  \align
  v&=\left(\tfrac{1}{2g-1}{v_g}^{2g-1}, \ldots, \tfrac13{v_g}^3, v_g\right), \\
  u\lr{j}&=\left(\tfrac1{2g-1}(u\lr{j}_g)^{2g-1}, \ldots, 
                  \tfrac13(u\lr{j}_g)^3, u\lr{j}_g\right). 
  \endalign
  $$
Then 
\fp
{\rm (1)} $S_{\sharp}(v)=-(-1)^{(g-1)(g-2)(g-3)/2}{v_g}^g$, \fp
{\rm (2)} $S_{\natural^{n+1}}(u\lr{1}+\cdots+u\lr{n}+v)$ \fp
$\qquad\qquad$  $=(-1)^{(g-n)(g-n-1)/2}S_{\natural^n}(u\lr{1}+\cdots+u\lr{n}){v_g}^{g-n}
            +(d^{\circ}(v_g)\geqq g-n+2)$, \fp
{\rm (3)} $S_{\flat}(2v)=-(-1)^{g(g-1)(g-2)/2}2{v_g}^{2g-1}$. 
\endproclaim

\demo{\it Proof} Since the statement (1) is proved by the same argument 
of our proof below of (2), we omit the proof.  
We prove (2) by separating into four cases with respect to 
the parity of  $g$  and  $n$.  
These four cases are similarly proved.  
Suppose both  $g$  and  $n$  are odd for instance.  
Then  $\natural^{n+1}$  is
  $$
  \{n+2, \ n+4, \ \ldots, \ g-2, \ g \}.
  $$
The number of the elements in  $\natural^{n+1}$  is  $(g-n)/2$.  
We denote this by  $\nu$ :  $\nu=(g-n)/2$.  
We are going to operating 
  $$
  D=\prod_{i\in \natural^{n+1}}\frac{\partial}{\partial u_i}
  =
  (2g-2n-1)\frac{\partial}{\partial p_{2g-2n+1}}
  (2g-2n-5)\frac{\partial}{\partial p_{2g-2n-3}}
  \cdots
          7\frac{\partial}{\partial p_7}
          3\frac{\partial}{\partial p_3}
  $$
on the determinant expression of  $S(u)$  
with respect to  $U_0$, $\ldots$, $U_{2g-1}$  
by using the formula in 2.3.  
The result of the operation is a sum of determinants 
with some of their rows are shifted several times to the right.  
We claim that {\it the terms in the sum except only one term, say  $|M|$,  
are vanish. }
Our proof of this claim is devided into three steps as follows.  

\fp{\it Step} 1. 
We write simply   $U_j=U\br{n+1}_j(\bk{u}_g, v_g)$.  
A determinant given by shifting one of the rows of determinant expression of  $S(u)$   
by two or more times by factors in  $D$  vanishes.  
Since this will be easily seen if the following argument is completed,  
we consider any determiant obtained by such a operation that the  $\nu$  factors in  $D$  above shift 
different  $\nu$  rows of  $S(u)$.  
Hence the factors of  $D$  changes  $\nu$  different  $U_{2j-1}$'s  to  $U_0$.  
We say that the entries in the first column to the  $(g-2\nu)^{\text{th}}$  column 
are in {\it the left region} of the matrix whose determinant expresses  $S(u)$  or of  $M$,   
and the other entries are in {\it the right region} of the matrix.  
The changed  $U_0$'s  must be exist in the right region.  
The reason is as follws.  
If any of the newly appeared  $\nu$ \  $U_0$'s are appeared in the left rigion
there exists in the right region at least one column 
which does not contain  $U_0$, 
and then the submatrix  $N_1$  consisting of the first column to it  
is a fundamental matrix without a simple row.  
Hence such a determiant vanishes by 1.3.  

\fp{\it Step} 2. We observe how a factor of  $D$  operates on  $|M|$.  
We pay attention to the  $(g-2\nu+1)^{\text{st}}$  column.  
We claim that the entry  $U_3$  in this column must be changed to  $U_0$  
by the factor  $\partial/\partial u_g=3\partial/\partial p_3$  of  $D$.  
If another entry  $U_{2j-1}(\neq U_3)$  in the column was changed to  $U_0$,  
of course by  $(2j-1)\partial/\partial p_{2j-1}$, 
we must use  $3\partial/\partial p_3$  to change  $U_3$  existing 
in another column.  
Then the entry  $U_0$  in just the row that contains the later  $U_3$ 
will be removed.  
Because this entry  $U_0$  is an entry in the right region, 
such an operation gives a determinant which contains 
a fundamental submatrix without a simple row 
with respect to  $u\lr{1}_g$, $\ldots$, $u\lr{n}_g$, $v_g$   
and it can not give  $|M|$  by 1.3.  
  $$
  \matrix
  \vdots & \vdots      & \vdots & \vdots & \vdots & \vdots &        \\
  U_4    & U_5         & U_6    & U_7    & U_8    & U_9    & \cdots \\
  U_2    & \boxed{U_3} & U_4    & U_5    & U_6    & U_7    & \cdots \\
  U_0    & U_1         & U_2    & U_3    & U_4    & U_5    & \cdots \\
         &             & U_0    & U_1    & U_2    & U_3    & \cdots \\
         &             &        &        & U_0    & U_1    & \cdots \\
         &             &        &        &        &        & \ddots \\
  \endmatrix
  $$
Now we look at the determinant  $3(\partial/\partial p_3)|M|$  
obtained above,  
and pay attention to the  $(2g-2\nu+3)^{\text{rd}}$  column of it.  
If  $7(\partial/\partial p_7)$  is used to this determinant 
without changing  $U_7$  in the row  to  $U_0$, 
this will change the other  $U_7$  in another row.  
Then the entry  $U_0$  in just the row where the later  $U_7$  exists 
will be removed and such operation of  $D$  gives a determinant 
which contains a fundamental submatrix without a simple row.  
Hence such operation gives a vanished determinant.  
  $$
  \matrix
  \vdots & \vdots & \vdots & \vdots      & \vdots & \vdots & \vdots & \vdots & \vdots &       \\
  U_4    & U_5         & U_6    & \boxed{U_7} & U_8    & U_9    &  *     &  *     &  *     & \cdots \\
         & \boxed{U_0} & U_1    & U_2         & U_3    & U_4    & U_5    & U_6    & U_7    & \cdots \\
  U_0    & U_1         & U_2    & U_3         & U_4    & U_5    & U_6    & U_7    &  *     & \cdots \\
         &             & U_0    & U_1         & U_2    & U_3    &  *     &  *     &  *     & \cdots \\
         &             &        &             & U_0    & U_1    &  *     &  *     &  *     & \cdots \\
         &             &        &             &        &        & U_0    & U_1    &  *     & \cdots \\
         &             &        &             &        &        &        &        & U_0    & \ddots \\
  \endmatrix
  $$
\fp{\it Step} 3. 
By the observation in the Step 2 
we see that there exists at most only one non-trivial determinant  $|M|$  in the sum,  
and that  $|M|$  is obtained by shifting rows of the determinant  $S(u)$  
such as just the  $(g, n+2)$,  $(g-2, n+1)$, $(g-3, n-1)$, $\ldots$, $(g-\nu+1, n-\nu+3)^{\text{rd}}$  
entries are changed to  $U_0$.  
The other three casese with respect to the parity of  $g$  and  $n$  can be 
similarly treated.  
Here we have completed the proof of the claim.  

From now on we consider all the four cases.  
Summing up the results above with being careful in the sign of the determinant, we have 
  $$
  \aligned
  S_{\natural^{n+1}}&(u_1+\cdots+u_n+v) \\
  &=\pm
   \left|\matrix
   U_g        & U_{g+1}    & \cdots & U_{g+n-1} & U_{g+n}   &  *     &  *     & \cdots & *      \\
   U_{g-2}    & U_{g-1}    & \cdots & U_{g+n-2} & U_{g+n-2} &  *     &  *     & \cdots & *      \\
   \vdots     & \vdots     & \ddots & \vdots    & \vdots    & \vdots & \vdots & \ddots & \vdots \\ 
   U_{g-2n-2} & U_{g-2n-1} & \cdots & U_{g-n+1} & U_{g-n-2} &  *     &  *     & \cdots & *      \\
   U_{g-2n}   & U_{g-2n+1} & \cdots & U_{g-n-1} & U_{g-n}   &  *     &  *     & \cdots & *      \\
              &            &        &           &           & U_0    &  *     & \cdots & *      \\
              &            &        &           &           &        & U_0    & \cdots & *      \\
              &            &        &           &           &        &        & \ddots & \vdots \\
              &            &        &           &           &        &        &        & U_0    \\
  \endmatrix\right|, 
  \endaligned
  \tag 2.1
  $$
where the sign  $\pm$  in the right hand side is 
  $$
  \alignat 2
   & (-1)^{(g-n+1)/2}   &   \ \ \ &\text{if} \ \  g-n\equiv 1  \ \text{mod} \ 2, \\ 
   &   1                &   \ \ \ &\text{if} \ \  g-n\equiv 0  \ \text{mod} \ 2.
  \endalignat
  $$
For the set  $\natural^n$  we have similarly
  $$
  S_{\natural^{n}}(u_1+\cdots+u_n) 
  =\pm 
  \left|\matrix
   U'_g        & U'_{g+1}    & \cdots & U'_{g+n-1} &  *     &  *     & \cdots & *      \\
   U'_{g-2}    & U'_{g-1}    & \cdots & U'_{g+n-2} &  *     &  *     & \cdots & *      \\
   \vdots      & \vdots      & \ddots & \vdots     & \vdots & \vdots & \ddots & \vdots \\ 
   U'_{g-2n-2} & U'_{g-2n-1} & \cdots & U'_{g-n-3} &  *     &  *     & \cdots & *      \\
               &             &        &            & U_0    &  *     & \cdots & *      \\
               &             &        &            &        & U_0    & \cdots & *      \\
               &             &        &            &        &        & \ddots & \vdots \\
               &             &        &            &        &        &        & U_0    \\
  \endmatrix\right|, 
  \tag 2.2
  $$
where  $U'_j=U\br{n+1}_j(\bk{u}_g, 0)=U\br{n}_j(\bk{u}_g)$, i.e. the values at  $v_g=0$, and 
the sign  $\pm$  is   
  $$
  \alignat 2
   &(-1)^{(g-n-1)/2}   &   \ \ \ &\text{if} \ \  g-n\equiv 1  \ \text{mod} \ 2,  \\
   &  1                &   \ \ \ &\text{if} \ \  g-n\equiv 0  \ \text{mod} \ 2.
  \endalignat
  $$
By writting the entries as polynomials of  $v_g$  and 
by using 1.3 to the coefficients of each power of  $v_g$  
we see that the determinant in the right hand side of  (2.1) is equal to
  $$
  \left|\matrix
   U'_g      & U'_{g+1}    & \cdots & U'_{g+n-1} & (-1)^{g+n}{v_g}^{g+n}     &  *     &  *     & \cdots & *      \\
   U'_{g-2}  & U'_{g-1}    & \cdots & U'_{g+n-3} & (-1)^{g+n-2}{v_g}^{g+n-2} &  *     &  *     & \cdots & *      \\
   \vdots    & \vdots      & \ddots & \vdots     & \vdots                    & \vdots & \vdots & \ddots & \vdots \\ 
   U'_{g-2n} & U'_{g-2n+1} & \cdots & U'_{g-n-1} & (-1)^{g-n}{v_g}^{g-n}     &  *     &  *     & \cdots & *      \\
             &             &        &            &                           & U_0    &  *     & \cdots & *      \\
             &             &        &            &                           &        & U_0    & \cdots & *      \\
             &             &        &            &                           &        &        & \ddots & \vdots \\
             &             &        &            &                           &        &        &        & U_0    \\
  \endmatrix\right|. 
  $$
The development along the $(n+1)^{\text{st}}$  column shows the desired formula of (2).  

If  $v=(\frac1{2g-1}{v_g}^{2g-1}, \ldots, \frac13{v_g}^3, v_g)$, 
then the formula in  (3) is given by a similar calculation in the proof of (2) as follows: 
  $$
  \align
  S_{\flat}(2v)
  &=(-1)^{g(g-1)(g-2)/2} \left(U\br{2}_g(v, v)\cdot U\br{2}_{g-1}(v, v)
   -U\br{2}_{g-2}(v, v)\cdot U\br{2}_{g+1}(v, v)\right) \\
  &=(-1)^{g(g-1)(g-2)/2}(-1)^{2g-1}\left((g+1){v_g}^g\cdot g{v_g}^{g-1} - (g-1){v_g}^{g-2}\cdot (g+2){v_g}^{g+1}\right) \\
  &=-(-1)^{g(g-1)(g-2)/2} \left((g+1)g - (g-1)(g+2)\right){v_g}^{2g-1} \\
  &=-(-1)^{g(g-1)(g-2)/2} 2{v_g}^{2g-1}.
  \endalign
  $$
Now the whole statements have been proved.  
\qedright
\enddemo

\newpage

\subheading{3. Hyperelliptic Functions}

\vskip 5pt
\fp
In this part we recall fundamentals 
of the theory of hyperelliptic functions. 

Let  $C$  be a smooth projective model of a curve of genus $g>0$ 
defined over  $\bold C$  whose affine equation is given 
by  $y^2=f(x)$, where
     $$
     f(x)=
     \lambda_0x^{2g+1} 
     +\lambda_1x^{2g}
     +\cdots
     +\lambda_{2g} x
     +\lambda_{2g+1}. 
     $$
In this paper, we keep the agreement  $\lambda_0=1$. 
We will use, however, the letter  $\lambda_0$  too when 
this notation makes easy to read an equation of homogeneous weight.  

We denote by  $\infty$  the point of  $C$  at infinity.  
It is known that the set of 
     $$
     \w_j
     :=\dfrac{x^{j-1}dx}{2y} \quad (j=1, \ldots, g) 
     $$
makes a basis of the space of the differential forms of the first kind.  
As usual we let  $\left[\omega' \ \omega''\right]$  be the period matrix 
for a suitable choice of the basis of the fundamental group of  $C$.  
Then the modulus of  $C$  is given by  
     $Z:={\omega'}^{-1}\omega''$.  
The lattice of periods is denoted by  $\Lambda$, that is 
     $$
     \Lambda
     := \omega' \tp{\left[\matrix \bold Z & \bold Z & \cdots & \bold Z
                  \endmatrix \right]}
     +  \omega''\tp{\left[\matrix \bold Z & \bold Z & \cdots & \bold Z
                  \endmatrix \right]}\ (\subset \bold C^g).
     $$
Let 
     $$
     \eta_j
     :=\dfrac{1}{2y}\sum_{k=j}^{2g-j}(k+1-j) \lambda_{2g-k-j} x^k dx
     \quad (j=1, \ldots, g), 
     $$
which are differential forms of the second kind without poles 
except at  $\infty$ 
(see  \cite{{\bf 2}}, p.195, Ex. i or \cite{{\bf 3}}, p.314). 
We introduce 
the matrices of periods  $\left[\eta' \  \eta''\right]$  
with respect to  $\eta_1$,  $\ldots$,  $\eta_g$   
for the basis of the fundamental group of  $C$  chosed in the above.  
We let 
     $$
     \delta'':
     =\ ^t\negthinspace
     \left[\matrix 
     \dfrac 12 &    \dfrac 12 & \cdots & \dfrac 12
     \endmatrix \right], \quad
     \delta':
     ={}^t\negthinspace
     \left[\matrix
     \dfrac g2 & \dfrac{g-1}2 & \cdots & \dfrac 12
     \endmatrix \right]
     \ \text{and}\ 
     \delta:
     =\left[\matrix 
     \delta''\\
     \delta'
     \endmatrix \right].
     $$
For  $a$  and  $b$  in  $\left(\frac 12 \bold Z\right)^g$, 
we let 
     $$
     \aligned
     &\vartheta\negthinspace
     \left[\matrix a \\ b \endmatrix\right](z)
     =\vartheta\negthinspace
     \left[\matrix a \\ b \endmatrix\right](z; Z) \\
     &=\sum_{n \in \bold Z^g} \exp \left[2\pi i\left\{
     \dfrac 12 \ ^t\negthinspace (n+a)Z(n+a) 
     + \ ^t\negthinspace (n+a)(z+b)\right\}\right].
     \endaligned
     $$
Then the hyperelliptic sigma function on  $\bold C^g$  
associated with  $C$  is defined by 
     $$
      \widetilde\sigma(u)
     =\ \text{exp}(-\dfrac{1}{2}u\eta'{\omega'}^{-1}\ ^t\negthinspace u)
     \vartheta\negthinspace
     \left[\delta\right]({\omega'}^{-1}\ ^t\negthinspace u;\ Z)
     $$
up to a multiplicative constant, 
where  $u=(u_1, u_2, \cdots, u_g)$.  We fix the constant latter. 

\definition{Definition 3.1} 
We define the {\it Sato weight}  by taking  
such the weight of  $u_j$  being  $2(g-j)+1$.  
Moreover we define the Sato weight of  $\lambda_j$  to be  $-2j$.  
\enddefinition

To fix the multiplicative constant above, we recall the following.  

\proclaim{Lemma 3.2}
{\rm (1)} The power-series development of  $\widetilde\sigma(u)$  
with respect to  $u_1$, $u_2$, $\cdots$, $u_g$  has the coefficients 
in polynomials of $\lambda_0$, $\lambda_1$, $\cdots$, $\lambda_{2g+1}$, 
and is homogeneous in Sato weight.  \fp
{\rm (2)} The terms of least total degree of 
the Taylor development of the function  $\widetilde\sigma(u)$  
at  $u=(0, 0, \ldots, 0)$  with respect to the variables  $u_1$, $\ldots$, $u_g$  
is a non-zero constant multiple of the Hankel type determinant
   $$
   \left|\matrix
   u_1         & u_2         & \cdots & u_{(g+1)/2} \\
   u_2         & u_3         & \cdots & u_{(g+3)/2} \\
   \vdots      & \vdots      & \ddots & \vdots      \\
   u_{(g+1)/2} & u_{(g+3)/2} & \cdots & u_g    
   \endmatrix\right|
   $$
if  $g$  is odd, or 
   $$
   \left|\matrix
   u_1         & u_2         & \cdots & u_{g/2} \\
   u_2         & u_3         & \cdots & u_{(g+2)/2} \\
   \vdots      & \vdots      & \ddots & \vdots      \\
   u_{g/2}     & u_{(g+2)/2} & \cdots & u_{g-1}    
   \endmatrix\right|
   $$
if  $g$  is even. 
\endproclaim

\demo{\it Proof}  
The statement (1) is just the Corollary 1 in \cite{{\bf 7}}, and 
(2) is proved in  \cite{{\bf 6}} p.32, Proposition 2.2 or \cite{{\bf 3}}, pp.359-360.   
\qedright
\enddemo

In this paper we let 
    $$
    \sigma(u)
    $$
be the function such that it is a constant multiple of  $\widetilde\sigma(u)$  
and  the terms of least total degree of 
its power-series development at  $u=(0, 0, \ldots, 0)$  
is just the Hankel type determinant above.  
It is easy to see from the proof of the Corollary 1 of \cite{{\bf 7}} 
that the power-series expansion of  $\sigma(u)$  
with respect to  $u_1$, $u_2$, $\cdots$, $u_g$  belongs to  
$\bold Q[\lambda_0, \lambda_1, \cdots, \lambda_{2g+1}][[u_1, u_2, \cdots, u_g]]$.

For  $u\in \bold C^g$  we conventionally denote by  $u'$  and  $u''$  
such elements of  $\bold R^g$  that 
  $u=\w'u'+\w''u''$,  
where  $\w'$  and  $\w''$  are those defined in Section 1. 
We define a  $\bold C$-valued  $\bold R$-bilinear form  $L(\quad, \quad)$  by
  $$
  L(u,v)=\tp{u}(\eta'v'+\eta''v'')
  $$
for  $u$, $v \in \bold C^g$.  
For  $\ell$  in  $\Lambda$, the lattice of periods as defined in Section 1, 
let 
  $$
  \chi(\ell)=\exp[2\pi i(\tp{\ell'}\delta''-\tp{\ell''}\delta')
  -\pi i\tp{\ell'}\ell''].
  $$

  \proclaim{Lemma 3.3}
  The function  $\sigma(u)$  is an odd function if  $g\equiv 1$  or  $2$  
  modulo  $4$,  
  and an even function if  $g\equiv 3$  or  $0$  modulo  $4$.  
  \endproclaim

\demo{\it Proof} 
See \cite{{\bf 16}} p.3.97 and p.3.100. 
\qedright
\enddemo

  \proclaim{Lemma 3.4}{\rm (the translational formula)}
  The function  $\sigma(u)$  satisfies 
  $$
  \sigma(u+\ell)=\chi(\ell)\sigma(u)\exp L(u+\tfrac 12\ell,\ell)
  $$
  for  all  $u\in \bold C^g$  and  $\ell\in \Lambda$. 
  \endproclaim

\fp
For a proof of this formula we refer to the reader to 
\cite{{\bf 2}}, p.286. 

\remark{\bf Remark 3.5} 
The Riemann form of  $\sigma(u)$  which is defined by 
$E(u,v)=L(u,v)-L(v,u)$,   ($u$, $v\in\bold C^g$) has simply 
written as   $E(u,v)=2\pi i(\tp{u'}v''-\tp{u''}v')$ 
(see \cite{{\bf 17}}, p.396, Lemma 3.1.2(2)).  
Hence, $E(\quad,\quad)$  is an  $i\bold R$-valued form 
and   $2\pi i \bold Z$-valued on  $\Lambda\times\Lambda$.  
Especially the pfaffian of  $E(\quad ,\quad )$  is  $1$.  
This simple expression is one of the good properties 
to distinguish  $\sigma(u)$  from the theta function 
without tuning by an exponential factor.  
\endremark

\ 

Let  $J$  be the Jacobian variety of the curve  $C$. 
We identify  $J$  with the Picard group  $\text{Pic}^{\circ}(C)$  of 
the linearly equivalent classes of divisors of degree zero of  $C$. 
Let  $\text{Sym}^g(C)$  be the  $g^{\text{th}}$  symmetric product of  $C$. 
Then we have a birational map 
  $$
  \align
  \text{Sym}^g(C) &\rightarrow \text{Pic}^{\circ}(C)=J \\
  (P_1, \ldots, P_g)&\mapsto \text{the class of}\ P_1+\cdots +P_g-g\cdot\infty.
  \endalign 
  $$ 
As an analytic manifold,  $J$  is identified with 
$\bold C^g/\Lambda$. We denote by  $\kappa$  the canonical map 
$\bold C^g \rightarrow \bold C^g/\Lambda=J$.  
We embed  $C$  into  $J$  by  
  $\iota$ : $Q \mapsto Q-\infty$.  
For each  $n=1$, $\ldots$, $g-1$  
let  $\Theta\br{n}$  be the subvariety of  $J$  
determined by the set of classes 
of the form  $P_1+\cdots +P_{n}-n\cdot \infty$  and is called 
the {\it standard theta subvariety} of dimension  $n$. 
Obviously  $\Lambda=\kappa^{-1}\iota(\infty)$  and  $\Theta\br{1}=\iota(C)$.  

Analytically each point  $(P_1, \ldots, P_g)$  of  $\text{Sym}^g(C)$  is 
canonically mapped to 
  $$
  u=(u_1, \ldots, u_g)=
  \left(\int_{\infty}^{P_1}+\cdots+\int_{\infty}^{P_g}\right)
  (\omega_1, \ \ldots,  \ \omega_g),  
  $$  
and  $\sigma(u)$  is regarded to be a function on 
the universal covering space  $\bold C^g$  of  $J$  with 
the canonical map  $\kappa$  above and the natural coordinate  $u$  of  $\bold  C^g$.  

If  $u=(u_1, \ \ldots, \ u_g)$  
is in  $\kappa^{-1}\iota(C)$,  we denote by  
  $$
  (x(u), y(u))
  $$  
the coordinate of the corresponding point on  $C$  by 
  $$
  u=
  \int_{\infty}^{(x(u), y(u))}
  (\frac{1}{2y}, \ \frac{x}{2y}, \ \ldots, \ \frac{x^{g-1}}{2y})dx 
  $$
with appropriate choice of a path of the integrals.  
Then we have  $(x(0, 0, \ldots, 0), \allowmathbreak y(0, 0, \ldots, 0))=\infty$.   
We frequently use the following lemma in the rest of the paper.  

  \proclaim{Lemma 3.6} Suppose  $u\in\kappa^{-1}\iota(C)$.  
  The Laurent expansion of  $x(u)$  and  $y(u)$  at  $u=(0, \ldots, 0)$  on the 
  pull-back  $\kappa^{-1}\iota(C)$  of  $C$  to  $\bold C^g$  are
  $$
  x(u)=\frac 1{u_g^2}+\Dg{u_g}{0}, \quad
  y(u)=-\frac 1{u_g^{2g+1}}+\Dg{u_g}{-2g+1}
  $$
  with their coefficients in  $\bold Q[\lambda_0, \lambda_1, \cdots, \lambda_{2g+1}]$.  
  Moreover  $x(u)$  and  $y(u)$  are homogeneous 
  of Sato weight  $-2$  and  $-(2g+1)$, respectively.     
  \endproclaim

\demo{\it Proof} 
We take  $t=\dfrac 1{\sqrt{x}}$  as a local parameter 
at  $\infty$  along  $\iota(C)$. 
If  $u$  is in  $\kappa^{-1}\iota(C)$  and 
sufficiently near  $(0, 0, \ldots, 0)$,  
we are agree to that  $t$, $u=(u_1, \ldots,  u_g)$  and  $(x, y)$  
are coordinates of the same point on  $C$. Then 
  $$
  \aligned
  u_g
  =&\int_{\infty}^{(x,y)}\frac {x^{g-1}dx}{2y} \\
  =&\int_{\infty}^{(x,y)}
   \frac{   x^{-3/2}dx   }
        {   
           2\sqrt{  
                     1+\lambda_1\frac 1x+\cdots
                     +\lambda_{2g+1}\frac 1{x^{2g+1}}
                  }
        } \\
  =&\int_0^t \frac{t^3\cdot\left(-\frac{2}{t^3}\right)dt}{2+\dg1} \\
  =&-t+\Dg{t}2.
  \endaligned
  $$
Hence  $x(u)=\frac {1}{u_g^2}+\Dg{u_g}{-1}$  and our assertion is proved, 
because  $x(-u)=x(u)$  and  $y(-u)=-y(u)$.  
The rest of the statements are obvious from the calculation above.  
\qedright
\enddemo

\proclaim{Lemma 3.7}
If  $u=(u_1, u_2, \ldots, u_g)$  is a variable on 
$\kappa^{-1}(\Theta\br{1})$,  then
  $$
  \aligned
  u_1&=\tfrac1{2g-1} {u_g}^{2g-1} + (d^{\circ}(u_g)\geqq 2g), \\ 
  u_2&=\tfrac1{2g-3} {u_g}^{2g-3} + (d^{\circ}(u_g)\geqq 2g-2), \\
  &\cdots\cdots\cdots \\
  u_{g-1}&=\tfrac13 {u_g}^3 + (d^{\circ}(u_g)\geqq 4)
  \endaligned
  $$
with the coefficients in  $\bold Q[\lambda_0, \lambda_1, \cdots, \lambda_{2g+1}]$,  
and these developments are homogeneous with respect to Sato weght.  
\endproclaim

\demo{\it Proof} 
The assertions are easily obtained by similar calculations 
as in the proof of 3.6.  
\qedright
\enddemo

\remark{\bf Remark 3.8}
The equalities in  (1.1)  are canonical limit relations of 
the equalities in 3.7 when we bring 
all the coefficients  $\lambda_1$, $\ldots$, $\lambda_{2g+1}$  close to  $0$, 
because of the homogeneousness in 3.6.  
\endremark

\newpage

\subheading{4. The Schur-Weierstrass Polynomial and the Sigma Function}

\vskip 5pt
\fp
The Schur-Weierstrass polynomial  $S(u)$  closely relates with 
the sigma function as follows.  We denote 
  $$
  \align
  S_{i_1i_2\cdots i_n}(u)&=
  \frac{\partial}{\partial u_{i_1}}
  \frac{\partial}{\partial u_{i_2}}
  \cdots
  \frac{\partial}{\partial u_{i_n}}
  S(u), \\
  \sigma_{i_1i_2\cdots i_n}(u)&=
  \frac{\partial}{\partial u_{i_1}}
  \frac{\partial}{\partial u_{i_2}}
  \cdots
  \frac{\partial}{\partial u_{i_n}}
  \sigma(u).
  \endalign
  $$

\proclaim{Propostion 4.1}
The power-series development of the function  $\sigma(u)$  
at  $u=(0, 0, \ldots, 0)$  is of the form 
  $$
  \sigma(u)
  =(-1)^{g(g-1)(g-3)/2}S(u)+(d^{\circ}(\lambda_1, \lambda_2, \ldots, \lambda_{2g+1})\geqq 1).  
  $$
\endproclaim

\demo{\it Proof} 
If  $g$  is odd then  
  $$
  S_{13\cdots g}(u)
  =\left|\matrix
      &     &     &     &        &        &        &        & U_0    \\
      &     &     &     &        &        & U_0    & U_1    & U_2    \\
      &     &     &     &        & \cdots & \vdots & \vdots & \vdots \\
      &     &     & U_0 & \cdots & \cdots & *      & *      & *      \\
      & U_0 & U_1 & *   & \cdots & \cdots & *      & *      & *      \\
  U_0 & U_1 & U_2 & *   & \cdots & \cdots & *      & *      & *      \\
      &     & U_0 & U_1 & \cdots & \cdots & *      & *      & *      \\
      &     &     &     &        & \ddots & U_1    & \vdots & \vdots \\
      &     &     &     &        &        &        & U_0    & U_1    
  \endmatrix
  \right|, 
  $$
and  $\sigma_{13\cdots g}(0, 0, \ldots, 0)=1$  by  3.1  with the definition of  $\sigma(u)$.  
If  $g$  is even then  
  $$
  S_{13\cdots (g-1)}(u)
  =\left|\matrix
      &     &     &     &        &        &        & U_0    \\
      &     &     &     &        & U_0    & U_1    & U_2    \\
      &     &     &     & \cdots & \vdots & \vdots & \vdots \\
      &     & U_0 & U_1 & \cdots & \cdots & *      & *      \\
  U_0 & U_1 & U_2 & *   & \cdots & \cdots & *      & *      \\
      & U_0 & U_1 & *   & \cdots & \cdots & *      & *      \\
      &     &     & U_0 & \cdots & \cdots & *      & *      \\
      &     &     &     &        & \ddots & \vdots & \vdots \\
      &     &     &     &        &        & U_0    & U_1    \\
  \endmatrix
  \right|, 
  $$
and  $\sigma_{13\cdots (g-1)}(0, 0, \ldots, 0)$  is also  $1$.  
Hence our statement follows from just the main theorem 
of \cite{{\bf 7}}, Theorem 6.3.  
\qedright
\enddemo

%
%

\proclaim{Corollary 4.2} 
For any derivatives of  $\sigma(u)$  and  $S(u)$  are related 
as follows\ {\rm :}
  $$
  \sigma_{i_1 i_2 \cdots i_n}(u)
  =(-1)^{g(g-1)(g-3)/2}S_{i_1 i_2 \cdots i_n}(u)
  +(d^{\circ}(\lambda_1, \lambda_2, \ldots, \lambda_{2g+1})\geqq 1).    
  $$
This series is homogeneous in Sato weight. 
\endproclaim

\demo{Proof} 
If a (higher) partial derivation of  $\sigma(u)$  with respect to  $u_1$, $\cdots$, $u_g$  
does not eliminate given terms, then it keeps order of Sato weight of the terms.  
Hence our statement follws from 3.2(1).  
\qedright
\enddemo

\newpage

\subheading{5. Vanishing Structure of the Sigma Function and of Its Derivatives}

\vskip 5pt
\fp
We investigate vanishing structure of  $\sigma(u)$  and  
of its  $derivatives$  by using 
Riemann singularity theorem and by a calculation of Brill-Noether matrices.  
The following is fundamental for us.  

\proclaim{Proposition 5.1}{\rm (Riemann singularity theorem)}  
For a given  $u\in \bold C^g$, we denote a divisor on the curve  $C$  
which gives the point  $u$  modulo  $\Lambda$  
by  $P_1+\cdots+P_g-g\cdot\infty$.   Then 
  $$
  \dim \Gamma(C, \Cal O(P_1+\cdots+P_g))=r+1
  $$
if and only if both of the following hold {\rm :}\fp
{\rm (1)} $\sigma_{i_1 i_2 \ldots i_h}(u)=0$  for any  $h<r$  
and for any  $i_1$, $\ldots$,  $i_h \in \{1, 2, \ldots, g\}$.  \fp
{\rm (2)} There exists an $r$-tuple  $\{i_1, i_2, \ldots, i_r\}$  such 
that  $\sigma_{i_1 i_2 \cdots i_r}(u)\neq 0$. 
\endproclaim

\demo{\it Proof} 
See  \cite{{\bf 1}}, p. 226. 
\qedright
\enddemo

To compute dimension of the  $0^{\text{th}}$  cohomology group above we recall 
Brill-Noether matrix defined as follows.  
We fix the local parameter of every point of  $C$.  
To make clear the following argument we define 
the local parameter  $t$  at each point  $P$  by 
  $$
  t=\cases
  x-x(P) \  &\text{if  $P$  is an ordinary point}, \\
    y    \  &\text{if  $P$  is a branch point different from  $\infty$}, \\
  \tfrac{1}{\sqrt{x}}
         \  &\text{if  $P=\infty$}.
  \endcases
  $$
Here we call  $P$  a branch point if  $y(P)=0$  or  $\infty$,  and 
an ordinary point otherwise. 

We denote by  $\Omega^1$  the sheaf of the differential forms of the first kind.  
For a point  $P$  of  $C$, let  $t$  be 
the local parameter defined above. We denote by  $P_t$  the point 
of  $C$  such that the value of  $t$  at  $P_t$  is  $t$. Then we define 
for  $\mu\in\Gamma(C,\Omega^1)$ 
  $$
  \delta^{\ell}\mu(P)
  =\frac{d^{\ell+1}}{dt^{\ell+1}}\int_{\infty}^{P_t}\mu\Big|_{t=0}. 
  $$
Since  $\mu$  is a holomorphic form, $\delta^{\ell}\mu(P)$  takes finite 
value at every point  $P$.  
Let  $D:=\sum_{j=1}^kn_jP_j$  be an effective divisor.  
We define by  $B(D)$  the matrix 
with  $g$  rows and  $\deg D:=\sum n_j$  columns 
whose  $(i, n_1+\cdots+n_{j-1}+\ell)$-entry is  $\delta^{\ell}\w_i(P_j)$.  
Then our computation starts by the following.  

\proclaim{Proposition 5.2}
Let  $D$  be an effective divisor of  $C$, 
Then
  $$
  \dim \Gamma (C, \Cal O(D)) = \deg D + 1 - \rank B(D).  
  $$
\endproclaim

\demo{\it Proof}  
For  $\mu\in\Gamma(C,\Omega^1)$,  
we can find uniquely  $c_1$, $\ldots$, $c_g\in \bold C$  such that 
$\mu=c_1\w_1+\cdots+c_g\w_g$.  
In this situation, the three statements \fp  
(1)  $\mu\in\Gamma(C, \Omega^1(-D))$, \fp
(2) $\delta^{\ell}\mu(P_j)=0$  for all  $j$  and  $\ell$  with 
    $1\leqq j\leqq k$  and  $0\leqq \ell\leqq n_j$, and \fp
(3) $B(D)\left[\matrix c_1 \\ \vdots \\ c_g \endmatrix\right]
        =\left[\matrix  0  \\ \vdots \\  0  \endmatrix\right]$  \fp
are equivalent. So 
  $\dim\Gamma(C, \Omega^1(-D))=g-\rank B(D)$.  
The Riemann-Roch theorem states
  $$
  \dim\Gamma(C, \Cal O(D))=\deg D -g+1+\dim\Gamma(C, \Omega^1(-D)).
  $$
Hence
  $$
  \dim\Gamma(C,\Cal O(D))=\deg D+1-\rank B(D). 
  $$
\qedright
\enddemo

To compute  $\rank B(D)$ we need only the case that  $D$  is of the form
  $$
  D=P_1+P_2+\cdots+P_n+(g-n)\cdot\infty, \ \ 
  (P_i\neq P_j  \text{ for any } i\neq j, \text{ and } P_j\neq \infty).
  $$  
Then the matrix  $B(D)$  is given by 
  $$
  \left[\matrix
  \frac1{2y}(P_1)         & \frac{x}{2y}(P_1)       & \cdots                  & 
  \frac{x^{n-1}}{2y}(P_1) & \vrule                  & \frac{x^n}{2y}(P_1)     & 
  \cdots                  & \frac{x^{g-3}}{2y}(P_1) & \frac{x^{g-2}}{2y}(P_1) & 
  \frac{x^{g-1}}{2y}(P_1) \\
\noalign{\vskip -2pt}
  \vdots                  & \vdots                  & \ddots                  & 
  \vdots                  & \vrule                  & \vdots                  & 
  \ddots                  & \vdots                  & \vdots                  & 
  \vdots                  \\
\noalign{\vskip -2pt}
  \frac1{2y}(P_n)         & \frac{x}{2y}(P_n)       & \cdots                  & 
  \frac{x^{n-1}}{2y}(P_n) & \vrule height17pt       & \frac{x^n}{2y}(P_n)     & 
  \cdots                  & \frac{x^{g-3}}{2y}(P_n) & \frac{x^{g-2}}{2y}(P_n) & 
  \frac{x^{g-1}}{2y}(P_n) \\
\noalign{\vskip -2pt}
                          &                         &                         & 
                          & \vrule height10pt       &                         & 
                          &                         &                         & 
                          \\
\noalign{\hrule}
                          &                         &                         & 
                          & \vrule height15pt       &         0               & 
  \cdots                  &       0                 &         0               & 
           0              \\
\noalign{\vskip -2pt}
                          &                         &                         & 
                          & \vrule height11pt       &         0               & 
  \cdots                  &       0                 &         0               & 
           1              \\
\noalign{\vskip -2pt}
                          &                         &                         & 
                          & \vrule\strut            &         0               & 
  \cdots                  &       0                 &         0               & 
                          \\
\noalign{\vskip -2pt}
                          &                         &                         & 
                          & \vrule\strut            &         0               & 
  \cdots                  &       0                 &         1               & 
                         \\
\noalign{\vskip -2pt}
                          &                         &                         & 
                          & \vrule\strut            &         0               & 
  \cdots                  &       0                 &                         & 
                          \\
\noalign{\vskip -2pt}
                          &                         &                         & 
                          & \vrule\strut            &         0               & 
  \cdots                  &       1                 &                         & 
                          \\
\noalign{\vskip -2pt}
                          &                         &                         & 
                          & \vrule\strut            &    \vdots               & 
  \ddots                  &                         &                         & 
                          \\
  \endmatrix
  \right].
  $$
Here the right low block was calculated by similar way as in 3.7.  
Therefore the rank  of  $B(D)$  is  $n+(g-n-1)/2$  or  $n+(g-n)/2$  
according as  $g-n$  is odd or even.  
Summing up the consideration above, we have 
  $$
  \align
  \dim \Gamma(C, \Cal O(P_1+P_2+\cdots+P_n+(g-n)\infty)) 
  &=g+1-(n+\lfloor(g-n-1)/2\rfloor) \\
  &=\lfloor(g-n+1)/2\rfloor+1.  
  \endalign
  $$
Again we denote by  $u$  a point in  $\bold C^g$  corresponding 
to  $P_1+P_2+\cdots+P_n+(g-n)\infty$.   
Proposition 5.1 yields that if   $h\leqq (g-n+1)/2$  
namely  $n\leqq g-2h-1$  then
  $$
  \sigma_{i_1 i_2 \cdots i_h}(u)=0 
  $$
for all  $i_1$, $\ldots$, $i_h$  and  
  $$
  \sigma_{j_1 j_2 \cdots j_{\lfloor(g-n+1)/2\rfloor}}(u)\neq 0 
  $$
for some  $j_1$, $j_2$, $\ldots$, $j_{\lfloor (g-n+1)/2\rfloor}$.  
Therefore we have the following.  

\proclaim{Lemma 5.3}
Fix arbitrarily  $h$  such that  $0\leqq h \leqq (g-1)/2$.  
For any integer  $n$  satisfying  $0\leqq n \leqq g-2h-1$  
the function   $u\mapsto \sigma_{i_1 i_2 \cdots i_h}(u)$   
on  $\kappa^{-1}(\Theta\br{n})$  is identically zero.  
\endproclaim

\newpage

\subheading{6. Special Derivatives of the Sigma function}

\vskip 5pt
\fp
We will introduce some special derivatives of the sigma function.  
These are important to state our Frobenius-Stickelberger type formula.  

\definition{Definition 6.1} 
Let  $\natural^n$  be the set defined in 2.1.   
Then we define a derivative  $\sigma_{\natural^n}(u)$  of  $\sigma(u)$  by 
  $$
  \sigma_{\natural^n}(u)
  =\left(\prod_{i\in \natural^n}
  \frac{\partial}{\partial u_i}\right)\negthinspace
  \sigma(u)
  $$
Especially we denote
  $$
  \sigma_{\sharp}(u)=\sigma_{\natural^1}(u), \ \ 
  \sigma_{\flat}(u)=\sigma_{\natural^2}(u).
  $$
\enddefinition

\fp
These functions are given as in the following table. 

$ $

\centerline{
\vbox{\offinterlineskip 
\halign{
\strut\vrule# & \ \hfil # \hfil \ &&
      \vrule# & \ \hfil # \hfil \ \cr 
\noalign{\hrule}
&   genus   && $\sigma_{\sharp}$ && $\sigma_{\flat}$ && $\sigma_{\natural^3}$ && $\sigma_{\natural^4}$ && $\sigma_{\natural^5}$ && $\sigma_{\natural^6}$ && $\sigma_{\natural^7}$ && $\sigma_{\natural^8}$ && $\sigma_{\natural^9}$ && $\sigma_{\natural^{10}}$ && $\cdots$  \cr
\noalign{\hrule}
\noalign{\hrule}
&     1     && $\sigma$          && ---              && ---                   && ---                   && ---                   && ---                   && ---                   && ---                   && ---                   && ---                   && $\cdots$  \cr
\noalign{\hrule}
&     2     && $\sigma_2$        && $\sigma$         && ---                   && ---                   && ---                   && ---                   && ---                   && ---                   && ---                   && ---                   && $\cdots$  \cr
\noalign{\hrule}
&     3     && $\sigma_2$        && $\sigma_3$       && $\sigma$              && ---                   && ---                   && ---                   && ---                   && ---                   && ---                   && ---                   && $\cdots$  \cr
\noalign{\hrule}
&     4     && $\sigma_{24}$     && $\sigma_3$       && $\sigma_4$            && $\sigma$              && ---                   && ---                   && ---                   && ---                   && ---                   && ---                   && $\cdots$  \cr
\noalign{\hrule}
&     5     && $\sigma_{24}$     && $\sigma_{35}$    && $\sigma_4$            && $\sigma_5$            && $\sigma$              && ---                   && ---                   && ---                   && ---                   && ---                   && $\cdots$  \cr
\noalign{\hrule}                                                                                                                           
&     6     && $\sigma_{246}$    && $\sigma_{35}$    && $\sigma_{46}$         && $\sigma_5$            && $\sigma_6$            && $\sigma$              && ---                   && ---                   && ---                   && ---                   && $\cdots$  \cr
\noalign{\hrule}
&     7     && $\sigma_{246}$    && $\sigma_{357}$   && $\sigma_{46}$         && $\sigma_{57}$         && $\sigma_6$            && $\sigma_7$            && $\sigma$              && ---                   && ---                   && ---                   && $\cdots$  \cr
\noalign{\hrule}
&     8     && $\sigma_{2468}$   && $\sigma_{357}$   && $\sigma_{468}$        && $\sigma_{57}$         && $\sigma_{68}$         && $\sigma_7$            && $\sigma_8$            && $\sigma$              && ---                   && ---                   && $\cdots$  \cr
\noalign{\hrule}
&     9    && $\sigma_{2468}$    && $\sigma_{3579}$  && $\sigma_{468}$        && $\sigma_{579}$        && $\sigma_{68}$         && $\sigma_{79}$         && $\sigma_8$            && $\sigma_9$            && $\sigma$              && ---                   && $\cdots$  \cr
\noalign{\hrule}                
& $\vdots$ && $\vdots$           && $\vdots$         && $\vdots$              && $\vdots$              && $\vdots$              && $\vdots$              && $\vdots$              && $\vdots$              && $\vdots$              && $\vdots$              && $\ddots$  \cr
}}}

\centerline{6.2. Table of  $\sigma_{\natural^n}(u)$}

$ $

\fp
We are going to prepare tools for investigation on the zeroes of these derivatives.  

\proclaim{Lemma 6.3} 
{\rm (1)} Let  $\check{\natural}^n$  be a proper subset of  $\natural^n$, and let
  $$
  \sigma_{\check{\natural}^n}(u)
  =\left(\prod_{i\in \check{\natural}^n}\frac{\partial}{\partial u_i}\right)
  \negthinspace\sigma(u)
  $$
Then the function  $u\mapsto  \sigma_{\check{\natural}^n}(u)$  
on  $\kappa^{-1}(\Theta\br{n})$  is identically zero. \fp
{\rm (2)} $\sigma_{\natural^{n+1}}(u)=0$  if and only if  $u\in\kappa^{-1}(\Theta\br{n})$.  
\endproclaim

\demo{Proof}
Obvious from 5.3
\qedright
\enddemo

\proclaim{Lemma 6.4}
Let  $n$  be an integer such that  $1\leqq n\leqq g-1$.   
Assume  $u$  belongs to  $\kappa^{-1}(\Theta\br{n})$.  
Then we have 
  $$
  \sigma_{\natural^n}(u+\ell)
  =\chi(\ell)\sigma_{\natural^n}(u)\exp L(u+\tfrac12\ell, \ell)
  $$
for all  $\ell\in \Lambda$.  
\endproclaim

\demo{Proof}
It follows from 6.3(2) 
after derivating the formula in 3.4 with respect to  $\natural^n$.  
\qedright
\enddemo

\proclaim{Proposition 6.5}
Let  $n$  be an integer such that  $1\leqq n \leqq g-1$. 
The space spanned by the functions  $u\mapsto \varphi(u)$  on   $\kappa^{-1}(\Theta\br{n})$  
vanishing only on  $\kappa^{-1}(\Theta\br{n-1})$  and  satisfying the equation 
  $$
  \varphi(u+\ell)
  =\chi(\ell)\varphi(u)\exp L(u+\tfrac12\ell, \ell)
  $$
for all  $\ell\in \Lambda$ is one dimensional.  
\endproclaim

\demo{\it Proof} 
Since  $\Theta\br{n-1}$  is a prime divisor of the variety  $\Theta\br{n}$, 
the vanishing order of such a function  $\varphi(u)$  is well-defined  
by the translational formula for  $\varphi(u)$.  
Let  $\varphi_1(u)$  and  $\varphi_2(u)$  be non-trivial functions 
on  $\kappa^{-1}(\Theta\br{n})$  with the properties of the statement. 
We may assume the vanishing order of  $\varphi_2(u)$  on  $\kappa^{-1}(\Theta\br{n-1})$  is 
less than or equal to that of  $\varphi_1(u)$.  
Then the function  $\varphi_1/\varphi_2$  is holomorphic on  $\kappa^{-1}(\Theta\br{n})$.  
Classically, here is a situation where we might use 
a special case of Hartogs' analytic continuation theorem.   
On the other hand, we have 
  $$
  \frac{\varphi_1}{\varphi_2}(u+\ell)
  =\frac{\varphi_1}{\varphi_2}(u)
  \ \ \ \text{for all  $u\in\kappa^{-1}(\Theta\br{n})$  and  $\ell\in \Lambda$},
  $$
by the translational formula.  
Therefore  $\varphi_1/\varphi_2$  can be regarded as a holomorphic function 
on  $\Theta\br{n}$.  Hence this is a constant function by Liouville's theorem. 
\qedright
\enddemo

\proclaim{Proposition 6.6} 
Let  $v$  be a variable on  $\kappa^{-1}(\Theta\br{1})$.  
\fp
{\rm (1)} Then the function  $v\mapsto \sigma_{\sharp}(v)$  has 
a zero of order  $g$  at  $v=(0, 0, \ldots, 0)$  modulo  $\Lambda$  
and no zero elsewhere. This function has a development of the form
  $$
  \sigma_{\sharp}(v)=(-1)^{(g-2)(g-3)/2}{v_g}^g+(d^{\circ}(v_g)\geqq g+2). 
  $$
{\rm (2)} Let  $n$  be an integer such that  $1\leqq n \leqq g-1$. 
Suppose   $v$, $u\lr{1}$, $u\lr{2}$, $\ldots$, $u\lr{n}$  
belong to  $\kappa^{-1}(\Theta\br{1})$. 
If  $u\lr{1}+\cdots+u\lr{n}\not\in\kappa^{-1}(\Theta\br{n-1})$, 
then the function  $v\mapsto \sigma_{\natural^{n+1}}(u\lr{1}+\cdots+u\lr{n}+v)$  
have zeroes of order  $1$  at  $v=-u\lr{1}$,  $\ldots$,  $-u\lr{n}$  
and a zero of order  $g-n$  at  $v=(0, 0, \ldots, 0)$  modulo  $\Lambda$  
and no other zero elsewhere.  
This function has a development of the form
  $$
  \sigma_{\natural^{n+1}}(u\lr{1}+\cdots+u\lr{n}+v)
           =(-1)^{(g-n)(g-n-1)/2}\sigma_{\natural^n}(u\lr{1}+\cdots+u\lr{n}){v_g}^{g-n}
            +(d^{\circ}(v_g)\geqq g-n+1).
  $$ 
\fp
{\rm (3)} If  $u\notin\kappa^{-1}(\Theta\br{g-1})$  then  $\sigma(u)\neq 0$.  
\endproclaim

\demo{\it Proof} The statement (3) is well-known 
(see \cite{{\bf 16}}, p.3.80, Theorem 5.3).  
The usual argument by integration of the logarithm of 
  $$
  \aligned
  \sigma_{\natural^{n+1}}&(u\lr{1}+\cdots+u\lr{n}+v+\ell) \\
  &=\chi(\ell)\sigma_{\natural^{n+1}}(u\lr{1}+\cdots+u\lr{n}+v)
  \exp L(u\lr{1}+\cdots+u\lr{n}+v+\tfrac12\ell, \ell)
  \endaligned
  \tag 6.1
  $$
of 6.4 along the boundary of a polygon representation of 
the Riemann surface of  $C$  shows that 
the functions  $v\mapsto \sigma_{\sharp}(v)$  
and  $v\mapsto \sigma_{\natural^{n+1}}(u\lr{1}+\cdots+u\lr{n}+v)$  
above have exactly  $g$  zeroes modulo  $\Lambda$  
or identically vanish.  
These functions, however, do not vanish identically because of 2.4 and 4.2.  
The other statements of (1) are followed from 4.2 and 2.4(1). 
\fp
We will prove the rest statement of (2). 
Lemma 6.3(2) yields that the function has zeroes at  $v=u\lr{1}$, $\ldots$, $u\lr{n}$.  
For  $v\in\kappa^{-1}(\Theta\br{1})$   we may write 
  $$
  \sigma_{\natural^{n+1}}(u\lr{1}+\cdots+u\lr{n}+v)
  =\sum_{j=0}^{\infty}\varphi\lr{j}(u\lr{1}+\cdots+u\lr{n}){v_g}^j 
  $$
with undetermined functions  $\varphi\lr{j}(u\lr{1}+\cdots+u\lr{n})$.  
(Then we know immediately that the term  $\varphi\lr{0}(u\lr{1}+\cdots+u\lr{n})$  vanishes by 6.3(2).) 
We denote by  $m$  the Sato weight of  $\sigma_{\natural^{n+1}}(u\lr{1}+\cdots+u\lr{n}+v)$.  
We see  $m=ng-\frac12n(n+1)$.  
The Sato weight of  $\varphi\lr{j}(u\lr{1}+\cdots+u\lr{n})$  is  $m-j$.  
The non-zero function  $u\mapsto \varphi\lr{j}(u)$  
with the smallest index  $j$  also satisfies the translational relation  
on  $\kappa^{-1}(\Theta\br{n})$  because of the Leibniz rule for 
operating  $(d/dv_g)^j$  to the translational formula 
of  $\sigma_{\natural^{n+1}}(u+v)$.    
We fix such  $j$.  
Then  $\varphi\lr{j}(u)$  is a consatnat multiple of  $\sigma_{\natural^n}(u)$  
on  $\kappa^{-1}(\Theta\br{n})$.  
However, such the constant is a polynomial of  $\lambda_0$, $\cdots$, $\lambda_{2g+1}$
by 4.2.  
Comparing Sato weight of  $\sigma_{\natural^{n+1}}(u+v)$  and that of  $\sigma_{\natural^n}(u)$ 
(that is equal to  $(n-1)g-\frac12 (n-1)n=m-(g-n)$ ), we have identically 
  $$
  \varphi\lr{j}(u)=0  \ \ \ \text{for} \ j=0, \ \ldots, \ g-n-1
  $$
on  $\kappa^{-1}(\Theta\br{n})$. 
Namely we see  $j=g-n$.   
Because we already found at least  $n$  zeroes 
of  $v\mapsto \sigma_{\natural^{n+1}}(u\lr{1}+\cdots+u\lr{n}+v)$  
for  $v\neq 0$,  it is impossible 
that  $\varphi\lr{g-n}(u\lr{1}+\cdots+u\lr{n})=0$  identically.  
So we have
  $$
  \sigma_{\natural^{n+1}}(u+v)=\varphi\lr{g-n}(u){v_g}^{g-n}+(d^{\circ}(v_g)\geqq g)
  $$
for non-trivial function  $\varphi\lr{g-n}(u)$.  
Moreover, by 6.3(2), we see  $\varphi\lr{g-n}(u)\neq 0$  if  $u\not\in\kappa^{-1}(\Theta\br{n-1})$.   
Thus we see by 2.4 and 4.2 that for  $u\in\kappa^{-1}(\Theta\br{n})$ 
  $$
  \varphi\lr{g-n}(u)=(-1)^{(g-n)(g-n-1)/2}\sigma_{\natural^n}(u).
  $$
Now all the other statements are clear.  
\qedright
\enddemo

\proclaim{Lemma 6.7} 
Let  $u\in\kappa^{-1}(\Theta\br{1})$.  Then
  $$
  \sigma_{\flat}(2u)=(-1)^{g-1} 2{u_g}^{2g-1}+\left(d^{\circ}(u_g)\geqq 2g+1\right). 
  $$
\endproclaim

\fp
{\it Proof.} 
The statement follows from 2.4(3),  4.2 and 3.3.  
\qedright

\proclaim{Lemma 6.8} 
Let  $u\in\kappa^{-1}(\Theta\br{1})$.  Then
  $$
  \frac{\sigma_{\flat}(2u)}{\sigma_{\sharp}(u)^4}=(-1)^g 2y(u).
  $$
\endproclaim

\fp
{\it Proof.} 
Lemma 6.4 shows that the left hand side is periodic with respect to  $\Lambda$, 
and is odd function by 3.3.  
We have by 6.7 that 
  $$
  \align
  \frac{\sigma_{\flat}(2u)}{\sigma_{\sharp}(u)^4} 
  &=\frac{(-1)^{g-1} 2{u_g}^{2g-1} + (d^{\circ}(u_g)\geqq 2g+1)}
         {\left({u_g}^g+(d^{\circ}(u_g)\geqq g+2)\right)^4} \\
  &=(-1)^{g-1}\frac{2}{{u_g}^{2g+1}}+\cdots \\
  &=(-1)^g 2y(u).
  \endalign
  $$
\qedright

\newpage

\subheading{7. Frobenius-Stickerberger Type Formulae}

\vskip 5pt
\fp
The initial case of our Frobenius-Stickelberger type formulae is 
as follows.  

\proclaim{Lemma 7.1} 
Suppose  $u$  and  $v$  are in  $\kappa^{-1}(\Theta\br{1})$.  
We have 
  $$
  (-1)^{g+1}\frac{\sigma_{\flat}(u+v) \sigma_{\flat}(u-v)}
       {\sigma_{\sharp}(u)^2\sigma_{\sharp}(v)^2}
  =-x(u)+x(v)\left
  (=\left|\matrix 1 & x(u) \\
                  1 & x(v) 
    \endmatrix\right|\right).
  $$
\endproclaim

\demo{Proof} As a function of  $u$  (or  $v$), 
we see the left hand side is periodic with respect to  $\Lambda$  by 6.4.  
Moreover we see that 
the left hand side has only pole at  $u=(0,0,\ldots, 0)$  modulo  $\Lambda$  by 6.6(1).  
Proposition 6.6 also shows the Laurent development of the left hand side is of the form 
  $$
  \frac{(\sigma_{\sharp}(v){u_g}^{g-1}+\cdots)
        (\sigma_{\sharp}(-v){u_g}^{g-1}+\cdots)}
       {({u_g}^g+\cdots)^2 \sigma_{\sharp}(v)^2} 
  =(-1)^g\frac{1}{{u_g}^2}+\cdots 
  =(-1)^gx(u)+\cdots.
  $$
Here we have used the fact that  $\sigma_{\sharp}(-v)=(-1)^g\sigma_{\sharp}(v)$  
which followes from 3.3.    
Since the both sides has the same zeroes at  $u=v$  and  $u=-v$, 
the two sides coincide.  
\qedright
\enddemo

The general case of our Frobenius-Stickelberger type formula is 
obtained as follows.  

\proclaim{Theorem 7.2} Let  $n$  be a fixed integer. 
Suppose  $u\lr{1}$, $\ldots$, $u\lr{n}$  are variables 
on  $\kappa^{-1}(\Theta\br{1})$.  
\fp
{\rm (1)} If  $1\leqq n \leqq g-1$, then we have
  $$
  c_n \frac{\sigma_{\natural^{n}}(u\lr{1}+\cdots+u\lr{n})
           \prod_{i<j}\sigma_{\flat}(u\lr{i}-u\lr{j})}
       {\sigma_{\sharp}(u\lr{1})^n\cdots\sigma_{\sharp}(u\lr{n})^n} 
  =\left|
  \matrix
  1      & x(u\lr{1}) & x^2(u\lr{1}) & \cdots & x^{n-1}(u\lr{1}) \\
  1      & x(u\lr{2}) & x^2(u\lr{2}) & \cdots & x^{n-1}(u\lr{2}) \\
  \vdots & \vdots     & \vdots       & \ddots & \vdots           \\
  1      & x(u\lr{n}) & x^2(u\lr{n}) & \cdots & x^{n-1}(u\lr{n}) 
  \endmatrix
  \right|, 
  $$
where  $c_n=(-1)^{g+1+\tfrac12(n-1)(n-2)(n-3)}$.  
\fp
{\rm (2)} If  $n \geqq g$, then we have
  $$
  \aligned
  &c_n\frac{\sigma(u\lr{1}+\cdots+u\lr{n})
           \prod_{i<j}\sigma_{\flat}(u\lr{i}-u\lr{j})}
       {\sigma_{\sharp}(u\lr{1})^n\cdots\sigma_{\sharp}(u\lr{n})^n} \\
  &=
  \left|
  \matrix
    1          & x(u\lr{1})       & x^2(u\lr{1}) & \cdots & x^g(u\lr{1}) 
  & y(u\lr{1}) & x^{g+1}(u\lr{1}) & xy(u\lr{1})  & x^{g+2}(u\lr{1}) 
  & \cdots \\
    1          & x(u\lr{2})       & x^2(u\lr{2}) & \cdots & x^g(u\lr{2}) 
  & y(u\lr{2}) & x^{g+1}(u\lr{2}) & xy(u\lr{2})  & x^{g+2}(u\lr{2}) 
  & \cdots \\
    \vdots     & \vdots           & \vdots       & \ddots & \vdots         
  & \vdots     & \vdots           & \vdots       & \vdots 
  & \ddots \\        
    1          & x(u\lr{n})       & x^2(u\lr{n}) & \cdots & x^g(u\lr{n}) 
  & y(u\lr{n}) & x^{g+1}(u\lr{n}) & xy(u\lr{n})  & x^{g+2}(u\lr{n})
  & \cdots 
  \endmatrix
  \right|,
  \endaligned
  $$
where the size of matrix of the right hand side is  $n\times n$  and  $c_n$  is given by the table below.  

\vskip 5pt

{\eightpoint
\centerline{
\vbox{\offinterlineskip 
\halign{
\strut\vrule# & \ \hfill # \hfil \ &&
      \vrule# & \ \hfill # \hfil \ \cr 
\noalign{\hrule}
& $g\backslash n$ {\rm mod} $4$  && $1$ && $2$ && $3$ && $0$  & \cr
\noalign{\hrule}
\noalign{\hrule}
&    $1$ \hfill                 &&  $1$ &&  $1$ && $-1$ && $-1$  & \cr 
\noalign{\hrule}
&    $2$ \hfill                 && $-1$ && $-1$ && $-1$ && $-1$  & \cr 
\noalign{\hrule}
&    $3$ \hfill                 && $-1$ &&  $1$ &&  $1$ && $-1$  & \cr 
\noalign{\hrule}
&    $0$ \hfill                 && $-1$ &&  $1$ && $-1$ &&  $1$  & \cr
\noalign{\hrule}                
}}}
}
\endproclaim

\demo{Proof} 
We regard two sides in each of cases (1) and (2) to be functions of  $u\lr{n}$.  
By 6.4 the left hand side is periodic with respect to  $\Lambda$.  
Hence we may regard the two sides of each case to be functions on the curve  $C$.  
It is easy to check the two sides have the same divisor by using 6.6 
and the theorem of Abel-Jacobi(see \cite{{\bf 19}}).  
The coefficients of the lowest terms of the Laurent developments of 
the two sides coincide by the hypothesis of induction.  
\qedright
\enddemo

\newpage

\subheading{8. Kiepert Type Formulae}
\vskip 5pt
\fp
The function  $\sigma(u)$  directly relates with  $x(u)$  as follows.  

\proclaim{Lemma 8.1} 
Fix  $j$  with  $0\leqq j \leqq g$.  
Let  $u$  and  $v$  are on  $\kappa^{-1}(\Theta\br{1})$.  
Then we have
  $$
  \lim_{u\to v}\frac{\sigma_{\flat}(u-v)}{u_j-v_j}
  =\frac 1{x^{j-1}(v)}.
  $$
\endproclaim

\demo{\it Proof} 
Because of 7.1 we have
  $$
  \frac{x(u)-x(v)}{u_j-v_j}
  =(-1)^g \frac{\sigma_{\flat}(u+v)}{\sigma_{\sharp}(u)^2\sigma_{\sharp}(v)^2}
   \cdot
   \frac{\sigma_{\flat}(u-v)}{u_j-v_j}.
  $$
Now we bring  $u$  close to  $v$.  
Then the limit of the left hand side is  
  $$
  \lim_{u\to v}\frac{x(u)-x(v)}{u_j-v_j}
  =\frac{dx}{du_j}(v).
  $$
This is equal to  ${2y}/{x^{j-1}}(v)$  by the definition.  
The required formula follows from 6.8. 
\qedright
\enddemo

\definition{Definition 8.2}
For  $u\in\kappa^{-1}(\Theta\br{1})$  we denote by  $\psi_n(u)$  the function 
${\sigma_{\natural^n}(nu)}/{\sigma_{\sharp}(u)^{n^2}}$  if  $n<g$, 
and  ${\sigma(nu)}/{\sigma_{\sharp}(u)^{n^2}}$  if  $n\geqq g$.   
\enddefinition

\fp
This function  $\psi_n(u)$  has the following expression,  
which is quite natural generalization of 
the classical formula of Kiepert \cite{{\bf 13}}.

\proclaim{Theorem 8.3} {\rm (Kiepert type formula)} We fix  $j$  with  $1\leqq j \leqq g$.  
Let  $u\in \kappa^{-1}(\Theta\br{1})$  and  $n$  be a positive integer. 
\fp
{\rm (1)} If  $1\leqq n \leqq g$  then  $\psi_n(u)=(-1)^{g+1+\frac12(n-1)(n^2-2)}(2y(u))^{n(n-1)2}$,  
and if  $n=g+1$  then  $\psi_n(u)=(2y(u))^{g(g+1)/2}$  or  $-(2y(u))^{g(g+1)/2}$  
according to  $g\equiv$ $0$ modulo $4$ or not.  
\fp
{\rm (2)} If  $n\geqq g$  then we have
  $$
  \aligned
  &c'_n(1!2!\cdots (n-1)!)\psi_n(u)
  =x^{(j-1)n(n-1)/2}(u)\times \\
  &\left|\matrix 
     x'           & (x^2)'            & \cdots        
  & (x^g)'        &  y'               & (x^{g+1})'     
  & (yx)'         & (x^{g+2})'        & \cdots      \\
     x''          & (x^2)''           & \cdots        
  & (x^g)''       &  y''              & (x^{g+1})''    
  & (yx)''        & (x^{g+2})''       & \cdots      \\
       x'''       & (x^2)'''          & \cdots        
  & (x^3)'''      &  y'''             & (x^{g+1})'''   
  & (yx)'''       & (x^{g+2})'''      & \cdots      \\
    \vdots           & \vdots               & \vdots        
  & \vdots           & \vdots               & \vdots        
  & \vdots           & \vdots               & \ddots     \\
    x\lr{n-1}     & (x^2)\lr{n-1}     & \cdots        
  & (x^g)\lr{n-1}    & y\lr{n-1}         & (x^{g+1})\lr{n-1} 
  & (yx)\lr{n-1}  & (x^{g+2})\lr{n-1} & \cdots     \\
  \endmatrix\right|(u),
  \endaligned
  $$
where the size of the matrix is  $n-1$  by  $n-1$,  
the symbols  ${}'$, ${}''$, $\ldots$, ${}\lr{n-1}$  denote  
$\frac{d}{du_j}$, $\left(\frac{d}{du_j}\right)^2$, $\ldots$, 
$\left(\frac{d}{du_j}\right)^{n-1}$, respectively, and 
  $$
  c'_n=\cases
   (-1)^{n-1}      & \ \ \ \text{if} \ \  g\equiv 1  \ \text{\rm mod} \ 4, \\ 
  -(-1)^{n(n-1)/2} & \ \ \ \text{if} \ \  g\equiv 2  \ \text{\rm mod} \ 4, \\ 
    -1             & \ \ \ \text{if} \ \  g\equiv 3  \ \text{\rm mod} \ 4, \\
   (-1)^{n(n+1)/2} & \ \ \ \text{if} \ \  g\equiv 0  \ \text{\rm mod} \ 4. 
  \endcases
  $$
\endproclaim

\demo{Proof}  
If  $1\leqq n \leqq g+1$, 
the right hand side of  7.2  is a determinant of Vandermonde.  
Hence we have the statement (1)  by using  8.1.  
Here we note that  $(-1)^{n(n-1)/2}c_n=c'_n$  if  $1\leqq n \leqq g$.  
The statement (2)  is proved by the same argument as in \cite{{\bf 18}} with using 7.2  and  8.1. 
\qedright
\enddemo

\remark{\bf Remark 8.4}
The polynomials  $\psi_n(u)$  are good generalization 
of {\it division polynomials} of an elliptic curve, 
and are useful to find torsion points on the curve  $C$ 
in the Jacobian variety  $J$.  
Indeed, for  $n\geqq g$, $u\in\kappa^{-1}\iota(C)$  is 
an  $n$-torsion in  $J$  if and only if 
all of  $\psi_{n-g+1}(u)$, $\psi_{n-g+2}(u)$, $\ldots$, $\psi_n(u)$, $\ldots$, 
$\psi_{n+g-1}(u)$  are vanish.  
Detailed description of this fact is seen in \cite{{\bf 8}}.  
\endremark

Finally we mention the degree of the polynomials above.  

\proclaim{Proposition 8.5}
The number of roots with counting multiplicities 
of the equation  $\psi_n(u)=0$  is  $\frac12n(n-1)(2g+1)$  if  $1\leqq n \leqq g-1$  
and  $n^2g-\frac12g(g+1)$  if  $n\geqq g$.  
\endproclaim

\demo{\it Proof} 
The number of the roots is equal to the order of the pole at  $u=(0,\ \ldots, 0)$.  
If  $u\in \kappa^{-1}(\Theta\br{1})$, then  $u_g$  is 
a local parameter at  $u=(0,\ \ldots, 0)$  
because of 3.7.  
We calculate it for the case  $j=g$.  
We denote  $u_g$  by  $t$.  
If  $n\geqq g$,  the lowest term of the Laurent expansion of the determinant in 8.3 
at  $u=(0, 0, \ldots, 0)$  is exactly the same as that of 
{\eightpoint
  $$
  \left|
  \matrix
    -\frac 2{t^3}  
  & -\frac 4{t^5}  
  & \cdots  
  & -\frac{2g}{t^{2g+1}}  
  & \cdots  
  & -\frac {n+g-1}{t^{n+g}} \\
     \frac{4\cdot 3}{t^4}  
  &  \frac{6\cdot 5}{t^6}  
  & \cdots  
  &  \frac{(2g+2)(2g+1)}{t^{2g+2}}  
  & \cdots  
  &  \frac{(n+g+1)(n+g)}{t^{n+g+1}} \\
    \vdots  
  & \vdots 
  & \ddots 
  & \vdots 
  & \ddots 
  & \vdots  \\
  \noalign{\vskip 5pt}
    (-1)^{n-1}\frac{(n+1)\cdots 3}{t^{n+1}}  
  & (-1)^{n-1}\frac{(n+3)\cdots 5}{t^{n+3}}
  & \cdots   
  & (-1)^{n-1}\frac{(n+2g-1)\cdots (2g-1)}{t^{2g+n-1}}  
  & \cdots  
  & (-1)^{n-1}\frac{(2n+g-1)\cdots (n+g-1)}{t^{2n+g-2}}  \\
  \endmatrix
  \right|.
  $$
}
This is equal to 
  $$
  \left|
  \matrix
  \dfrac1{t^3} &   &       &               &        &           &          \\
      & \dfrac1{t^6}       &               &        &           &          \\
      &           & \ddots &               &        & \text{\Large  *} &   \\
      &           &       & \dfrac1{t^{3g}} &       &           &          \\
      &           &        &       &\dfrac1{t^{3g+2}} &         &          \\
      &           &        &               &        & \ddots    &          \\
      &           &        &               &        &  & \dfrac1{t^{2n+4g}} \\
  \endmatrix
  \right|
  $$
times a non-zero constant.  
This determinant is  $(\frac 32 g(g+1)+(n-g-1)(n+2g))^{\text{th}}$  power of  $1/t$.  
We see  $x^{(g-1)n(n-1)/2}(u)=1/t^{(g-1)n(n-1)}+\cdots$.  
So the lowest term of the Laurent expansion of the right hand side of 
the formula in 8.3 is  $\left(n^2g-\frac12 g(g+1)\right)^{\text{st}}$  power of  $1/t$.  
If  $n<g$  then the lowest term is easily seen to be  $\frac12 n(n-1)(2g+1)^{\text{st}}$  
power of  $1/t$.  
\qedright
\enddemo

\newpage

\leftheadtext\nofrills{\eightpoint Shigeki Matsutani}

\specialhead{\bf Appendix:  
Connection of The formulae of Cantor-Brioschi and of Kiepert type  \ 
(by S. Matsutani). } \endspecialhead

\vskip 5pt
\fp
In this Appendix we prove a formula of Cantor in \cite{{\bf 8}} 
(Theorem A.1 below) by using 8.3.  
This is detailed exposition of the appendix in \cite{{\bf 15}}.  
Since our argument is convertible, 
8.3 is proved by using the formula of Cantor.  

Let  $u=(u_1, u_2, \ldots, u_g)$  be the system of 
variables explained in Section 3.  
We assume that  $u$  belongs to  $\kappa^{-1}\iota(C)$.  
So we may use the notation  $x(u)$  and  $y(u)$.  
If  $\mu(u)$  is a function on  $\kappa^{-1}\iota(C)$  we can regard it 
locally as a function of  $u_1$.  
We denote by
  $$
  \mu'(u), \ 
  \mu''(u), \ 
  \ldots, \  
  \mu\lr{\nu}(u), \  
  \ldots
  $$
the functions obtained by operating  
  $$
  \tfrac d{du_1}, \ 
  (\tfrac d{du_1})^2, \
  \ldots, \  
  (\tfrac d{du_1})^{\nu}, \ 
  \ldots
  $$
for the function  $\mu(u)$  along  $\iota(C)$;   
and by
  $$
  \dot{\mu}(u), \ 
  \ddot{\mu}(u), \ 
  \ldots, \  
  \mu\do{\nu}(u), \  
  \ldots
  $$
the functions given by operating  
  $$
  \tfrac d{dx}, \ 
  (\tfrac d{dx})^2, \
  \ldots, \  
  (\tfrac d{dx})^{\nu}, \ 
  \ldots
  $$
for  $\mu(u)$.  
Here we regard  $\mu(u)$  locally as a function of  $x=x(u)$.  

\fp
Recall  $\psi_n(u)$  defined in 8.2. 
The determinant expression of  $\psi_n(u)$  due to Cantor is following.  

\proclaim{Theorem A.1} {\rm (Cantor \cite{{\bf 8}})} 
Suppose  $n\geqq g$.  
Let  $r$  be the largest integer not excced  $(n-g-1)/2$, 
and  $s=n-1-r$.  Then 
  $$
  \allowdisplaybreaks
  \psi_n(u) 
  =\varepsilon_n \cdot (2y)^{n(n-1)/2} 
  \times
  \cases
  \left|
  \matrix
    \dfrac {y\do{g+2}}{(g+2)!} 
  & \dfrac {y\do{g+3}}{(g+3)!}
  & \cdots 
  & \dfrac {y\do{r+1}}{(r+1)!} \\ 
    \dfrac {y\do{g+3}}{(g+3)!} 
  & \dfrac {y\do{g+4}}{(g+4)!}
  & \cdots 
  & \dfrac {y\do{r+2}}{(r+2)!} \\ 
    \vdots
  & \vdots
  & \ddots
  & \vdots \\
    \dfrac {y\do{r+1}}{(r+1)!} 
  & \dfrac {y\do{r+2}}{(r+2)!}
  & \cdots 
  & \dfrac {y\do{n-1}}{(n-1)!} 
  \endmatrix
  \right|,
  &  \text{if \ \  $n\not\equiv g \mod 2${\rm;}} \\
  & \strut \\
  \left|
  \matrix
    \dfrac {y\do{g+1}}{(g+1)!} 
  & \dfrac {y\do{g+2}}{(g+2)!}
  & \cdots 
  & \dfrac {y\do{r+1}}{(r+1)!} \\ 
    \dfrac {y\do{g+2}}{(g+2)!} 
  & \dfrac {y\do{g+3}}{(g+3)!}
  & \cdots 
  & \dfrac {y\do{r+2}}{(r+2)!} \\ 
    \vdots
  & \vdots
  & \ddots
  & \vdots \\
    \dfrac {y\do{r+1}}{(r+1)!} 
  & \dfrac {y\do{r+2}}{(r+2)!}
  & \cdots 
  & \dfrac {y\do{n-1}}{(n-1)!} 
  \endmatrix
  \right|,
  &  \text{if \ \   $n\equiv g \mod 2$, } 
  \endcases
  $$
where  $\varepsilon_n$  is given by the table below.

\vskip 5pt

\centerline{
\vbox{\offinterlineskip 
\halign{
\strut\vrule# & \ \hfill # \hfil \ &&
      \vrule# & \ \hfill # \hfil \ \cr 
\noalign{\hrule}
& $g$ {\rm mod} $4$ $\backslash$ $n$ {\rm mod} $8$ 
                                &&  $1$ &&  $2$ &&  $3$ &&  $4$ 
                                &&  $5$ &&  $6$ &&  $7$ &&  $0$  & \cr
\noalign{\hrule}
\noalign{\hrule}
&    $1$ \hfill                 &&  $1$ &&  $1$ && $-1$ &&  $1$ 
                                &&  $1$ &&  $1$ &&  $1$ && $-1$  & \cr 
\noalign{\hrule}
&    $2$ \hfill                 && $-1$ && $-1$ &&  $1$ &&  $1$ 
                                && $-1$ && $-1$ &&  $1$ &&  $1$  & \cr 
\noalign{\hrule}
&    $3$ \hfill                 && $-1$ &&  $1$ && $-1$ &&  $1$ 
                                && $-1$ &&  $1$ && $-1$ &&  $1$  & \cr 
\noalign{\hrule}
&    $0$ \hfill                 && $-1$ &&  $1$ && $-1$ && $-1$ 
                                && $-1$ &&  $1$ &&  $1$ &&  $1$  & \cr
\noalign{\hrule}                
}}}

\endproclaim

\remark{\bf Remark A.2} 
{\rm (1)} The both matrices above are of size  $s\times s$.  
\fp
{\rm (2)} The number  $r$  above is the largest number  $k$  of 
the entries of the form  $(x^k)'$  in the right hand side of the formula in 8.3. 
\fp
{\rm (3)} The constant factor is not clear 
by our definition in 8.2 and  \cite{{\bf 8}} of  $\psi_n(u)$.  
It is determined when our calculation has completed.  \fp
\fp
{\rm (4)} This formula for the case  $g=1$  is classically known by Brioschi \cite{{\bf 5}}.  
\endremark

\nopagebreak
$ $

The following Lemma is easily checked.  

\proclaim{Lemma A.3}
Let  $m>0$  be any integer.  One has
  $$
  \left(\frac d{du_1}\right)^m
  =(2y)^m(u)\left(\frac d{dx}\right)^m
  +\sum_{j=1}^m {a\lr{m}_j}(u)\left(\frac d{dx}\right)^j.
  $$
Here  $a\lr{m}_j(u)$  are polynomials 
of  $y(u)$, $\frac{dy}{dx}(u)$,  $\frac{d^2y}{dx^2}(u)$, $\ldots$, $\frac{d^{m-1}y}{dx^{m-1}}(u)$.  
\endproclaim

Let  $s=n-1-r$. Then,  by 8.3 for  $j=1$, we have
  $$
  \aligned
  &c'_n 1!2!\cdots (n-1)!\psi_n(u) \\
  & =(-1)^{r(r+1)/2}\left|\matrix 
     x'           & (x^2)'            & \cdots       & (x^r)'  
  &  y'           & (yx)'             & \cdots       & (yx^{s-1})'    \\
     x''          & (x^2)''           & \cdots       & (x^r)'' 
  &  y''          & (yx)''            & \cdots       & (yx^{s-1})''   \\
     x'''         & (x^2)'''          & \cdots       & (x^r)''' 
  &  y'''         & (yx)'''           & \cdots       & (yx^{s-1})'''  \\ 
    \vdots        & \vdots            &              & \vdots  
  & \vdots        & \vdots            & \ddots       & \vdots     \\
    x\lr{n-1}     & (x^2)\lr{n-1}     & \cdots       & (x^r)\lr{n-1}
  & y\lr{n-1}     & (yx)\lr{n-1}      & \cdots       & (yx^{s-1})\lr{n-1} 
  \endmatrix\right|(u).
  \endaligned
  $$
Here the size of the matrix is  $n-1$  by  $n-1$.  
By A.3 we have
  $$
  \left[
  \matrix 
  \frac d{du_1} \\
  \left(\frac d{du_1}\right)^2 \\
  \left(\frac d{du_1}\right)^3 \\
  \vdots \\
  \left(\frac d{du_1}\right)^{n-1} 
  \endmatrix
  \right]
  =
  \left[
  \matrix
        2y      &               &               &           &           \\
  a\lr{1}_2     & (2y)^2        &               &           &           \\
  a\lr{1}_3     & a\lr{2}_3     & (2y)^3        &           &           \\
     \vdots     &   \vdots      &    \vdots     & \ddots    &           \\
  a\lr{1}_{n-1} & a\lr{2}_{n-1} & a\lr{3}_{n-1} & \cdots    & (2y)^{n-1}
  \endmatrix
  \right]
  \left[
  \matrix
  \frac d{dx} \\
  \left(\frac d{dx}\right)^2 \\
  \left(\frac d{dx}\right)^3 \\
  \vdots \\
  \left(\frac d{dx}\right)^{n-1}
  \endmatrix
  \right].
  $$
Now we consider 
  $$
  \tp{
  \left[
  \matrix
  \frac d{dx}\mu &
  \left(\frac d{dx}\right)^2 \mu &
  \left(\frac d{dx}\right)^3 \mu &
  \cdots &
  \left(\frac d{dx}\right)^{n-1} \mu
  \endmatrix
  \right]}
  $$
for  $\mu=x$, $x^2$, $\ldots$  and  $y$, $yx$, $yx^2$, $\ldots$.  
Obviously
  $$
  \left[
  \matrix
  \frac d{dx} \\
  \left(\frac d{dx}\right)^2  \\
  \left(\frac d{dx}\right)^3  \\
  \vdots \\
  \left(\frac d{dx}\right)^{n-1} 
  \endmatrix
  \right]
  \left[
  \matrix
  x & x^2 & \cdots & x^r
  \endmatrix
  \right]
  =
  \left[
  \matrix
  1! &    &    &            &     \\
     & 2! &    & \text{\Large *} &     \\
     &    & 3! &            &     \\ 
     &    &    & \ddots     &     \\
     &    &    &            & r!  \\
     \noalign{\vskip-7pt}
     \multispan5\hrulefill        \\
     &    &    &            &     \\
     &    & \text{\Large 0} &   &  \\
     &    &    &            & 
  \endmatrix
  \right].
  $$
For  $\mu=$, $y$, $yx$, $\ldots$, $yx^{s-1}$, we have  
  $$
  \allowdisplaybreaks
  \align
  &
  \left[
  \matrix
  \frac d{dx} \\
  \left(\frac d{dx}\right)^2  \\
  \left(\frac d{dx}\right)^3  \\
  \vdots \\
  \left(\frac d{dx}\right)^{n-1} 
  \endmatrix
  \right]
  \left[
  \matrix
  y & yx & yx^2 & \cdots & yx^{s-1} 
  \endmatrix
  \right] \\
  {\eightpoint =}&
  {\eightpoint 
  \left[
  \matrix
    \C{1}{0} \dot{y} 
  & \C{1}{0} \dot{y}x  +\C{1}{1}y \hfill
  & \C{1}{0} \dot{y}x^2+\C{1}{1}y\cdot 2x \hfill  
  & \C{1}{0} \dot{y}x^3+\C{1}{1}y\cdot 3x^2 \hfill
  & \cdots                                   \\
    \C{2}{0} \ddot{y} 
  & \C{2}{0} \ddot{y}x  +\C{2}{1} \dot{y} 
  & \C{2}{0} \ddot{y}x^2+\C{2}{1} \dot{y}\cdot 2x  +\C{2}{2}y\cdot 2!  
  & \C{2}{0} \ddot{y}x^3+\C{2}{1} \dot{y}\cdot 3x^2
    +\C{2}{2}y\cdot 3\cdot 2x \hfill
  & \cdots                                   \\
    \C{3}{0}\dddot{y} 
  & \C{3}{0}\dddot{y}x  +\C{3}{1}\ddot{y} 
  & \C{3}{0}\dddot{y}x^2+\C{3}{1}\ddot{y}\cdot 2x  +\C{3}{2}\dot{y}\cdot 2!  
  & \C{3}{0}\dddot{y}x^3+\C{3}{1}\ddot{y}\cdot 3x^2
    +\C{3}{2}\dot{y}\cdot 3\cdot 2x +\C{3}{3}y\cdot 3! 
  & \cdots                                   \\
    \vdots
  & \vdots
  & \vdots
  & \vdots
  & \ddots                                   \\
  \endmatrix
  \right] 
} \\ 
  =&
  \left(
  \left[
  \matrix
    \C{1}{0}\Dxd       & \C{1}{1}           &                    
  &                    &                    \\
    \C{2}{0}\Dx{2}     & \C{2}{1}\Dxd       & \C{2}{2}           
  &                    &                    \\
    \vdots             & \vdots             & \vdots             
  & \ddots             &                    \\
    \C{s-1}{0}\Dx{s-1} & \C{s-1}{1}\Dx{s-2} & \C{s-1}{2}\Dx{s-3} 
  & \cdots             & \C{s-1}{s-1}       \\
    \vdots             & \vdots             & \vdots             
  & \ddots             & \vdots             \\
    \C{r}{0}\Dx{r}     & \C{r}{1}\Dx{r-1}   & \C{r}{2}\Dx{r-2}   
  & \cdots             & \C{r}{s-1}\Dx{r-s+1} \\
  \noalign{\vskip 0pt}  
  \multispan5\hrulefill        \\
    \C{r+1}{0}\Dx{r+1} & \C{r+1}{1}\Dx{r}   & \C{r+1}{2}\Dx{r-1} 
  & \cdots             & \C{r+1}{s-1}\Dx{r-s+2} \\
    \vdots             & \vdots             & \vdots             
  & \ddots             & \vdots             \\
    \C{n-1}{0}\Dx{n-1} & \C{n-1}{1}\Dx{n-2} & \C{n-1}{2}\Dx{n-3} 
  & \cdots             & \C{n-1}{s-1}\Dx{n-s} \\
  \endmatrix\right]
  \right.   \\
  & \times
  \left.
  \left[\matrix
    y & yT           & y\cdot T^2   & \cdots     & y\cdot T^{s-1}           \\
      & y            & y\cdot 2T    & \cdots     & y\cdot (s-1)T^{s-2}      \\
      &              & y\cdot 2!    & \cdots     & y\cdot (s-1)(s-2)T^{s-3} \\
      &              &              & \ddots     & \vdots                   \\
      &              &              &            & y\cdot (s-1)!            \\
  \endmatrix\right]\right)\Bigg|_{T=x}.
  \endalign
  $$
Thus we have
  $$
  \det
  \left(
  \left[
  \matrix
  \frac d{dx} \\
  \left(\frac d{dx}\right)^2  \\
  \left(\frac d{dx}\right)^3  \\
  \vdots \\
  \left(\frac d{dx}\right)^{n-1} 
  \endmatrix
  \right]
  \left[\matrix
   x & x^2 & \cdots & x^r & y & yx & \cdots & yx^{s-1} 
  \endmatrix\right]\right)
  $$
is equal to  $(1!2!\cdots r!)$  times
{\eightpoint
\par
  $$
  \allowdisplaybreaks
  \align
  &\det
  \left(
  \left[
  \matrix
  \C{r+1}{0}\Dx{r+1} & \cdots & \C{r+1}{s-1}\Dx{r-s+2} \\
  \vdots             & \ddots & \vdots                 \\
  \C{n-1}{0}\Dx{n-1} & \cdots & \C{n-1}{s-1}\Dx{n-s}   \\
  \endmatrix\right]
  \left[\matrix
  y & yT & y\cdot T^2 & \cdots & y\cdot T^{s-1}           \\
    & y  & y\cdot 2T  & \cdots & y\cdot (s-1)T^{s-2}      \\
    &    & y\cdot 2!  & \cdots & y\cdot (s-1)(s-2)T^{s-3} \\
    &    &            & \ddots & \vdots                   \\
    &    &            &        & y\cdot (s-1)!            \\
  \endmatrix\right]
  \right)\Bigg|_{T=x} \\
  =&
  \det
  \left(
  \left[
  \matrix
  \C{r+1}{0}\Dx{r+1} & \cdots & \C{r+1}{s-1}\Dx{r-s+2} \\
  \vdots             & \ddots & \vdots                 \\
  \C{n-1}{0}\Dx{n-1} & \cdots & \C{n-1}{s-1}\Dx{n-s}   \\
  \endmatrix\right]
  \left[\matrix
  y &           &            &        &                 \\
    & y\cdot 1! &            &        &                 \\
    &           & y\cdot 2!  &        &                 \\
    &           &            & \ddots &                 \\
    &           &            &        & y\cdot (s-1)!   \\
  \endmatrix\right]
  \right)  \\
  =&
  \det
  \left(
  \left[
  \matrix
    \C{r+1}{0}\Dx{r+1}   & 1!\C{r+1}{1}\Dx{r}    
  & \cdots & (s-1)!\C{r+1}{s-1}\Dx{r-s+2}  \\
    \C{r+2}{0}\Dx{r+2}   & 1!\C{r+2}{1}\Dx{r+1}    
  & \cdots & (s-1)!\C{r+2}{s-1}\Dx{r-s+3}  \\
    \vdots               & \vdots                
  & \ddots & \vdots                        \\
    \C{n-1}{0}\Dx{n-1}   & 1!\C{n-1}{1}\Dx{n-2}  
  & \cdots & (s-1)!\C{n-1}{s-1}\Dx{n-s}    \\
  \endmatrix\right]
  \left[\matrix
  y &   &   &        &    \\
    & y &   &        &    \\
    &   & y &        &    \\
    &   &   & \ddots &    \\
    &   &   &        & y  \\
  \endmatrix\right] 
  \right). 
  \endalign
  $$
} 
By dividing the first row by  $(r+1)!$, 
              the second row by $(r+2)!$, and so on,  
we see the above is equal to
  $$
  (r+1)!(r+2)!\cdots (n-1)!
  \left|
  \matrix
    \dfrac {y\do{r+1}}{(r+1)!} 
  & \dfrac {y\do{r}}{r!}
  & \cdots  &  \dfrac {y\do{r-s+2}}{(r-s+2)!} \\
    \dfrac {y\do{r+2}}{(r+2)!} 
  & \dfrac {y\do{r+1}}{(r+1)!}
  & \cdots  &  \dfrac {y\do{n-s+3}}{(r-s+3)!}   \\
    \vdots
  & \vdots
  & \ddots  &  \vdots                         \\
    \dfrac {y\do{n-1}}{(n-1)!} 
  & \dfrac {y\do{n-2}}{(n-2)!}
  & \cdots & \dfrac {y\do{n-s}}{(n-s)!}       \\
  \endmatrix
  \right|
  $$
Hence we obtained the determinant in A.1 by reshuffling the columns. 
Then the factor  $(-1)^{n-r-1}=(-1)^s$  appears.  
Summing up the calculation above and checking  
  $$
  c'_n\cdot (-1)^{s+\frac12r(r-1)}=\varepsilon_n,
  $$  
we obtained Theorem A.1.  

\leftheadtext\nofrills{\eightpoint Yoshihiro \^Onishi}

\vskip 10pt

{\eightpoint 
\fp
Yoshihiro \^Onishi
\fp
Faculty of Humanities and Social Sciences, Iwate University
\fp
Morioka, 020-8550, Japan 
\fp
{\it E-mail address:} onishi\@iwate-u.ac.jp

\vskip 8pt

\fp
Shigeki Matsutani
\fp
8-21-1, Higashi-Linkan Sagamihara, 228-0811, Japan
\fp
{\it E-mail address:} RXB01142\@nifty.or.jp
}

\enddocument
\bye